\documentclass[12pt]{amsart}
\usepackage[latin1]{inputenc}
\usepackage{amsmath}
\usepackage{amssymb}
\usepackage{graphicx}
\usepackage{graphics}
\usepackage{amsthm}
\usepackage{amscd}
\usepackage{color}
\usepackage{amsfonts}
\usepackage{fancyhdr}
\newtheorem{theorem}{Theorem}[section]
\newtheorem{mainthm}{Theorem}
\newtheorem*{theorem*}{Theorem}
\newtheorem{corollary}[theorem]{Corollary}
\newtheorem{proposition}[theorem]{Proposition}
\newtheorem{lemma}[theorem]{Lemma}
\newtheorem{Remark}[theorem]{Remark}
\newtheorem*{definition*}{Definition}

\newtheorem{claim}[theorem]{Claim}
\newtheorem{definition}[theorem]{Definition}
\newtheorem*{Question*}{Question}
\newtheorem{Question}{Question}

\def\N{\mathbb{N}}
\def\Z{\mathbb{Z}}

\def\R{\mathbb{R}}

\def\T{\mathbb{T}}

\def\norm #1{\Vert \,#1\, \Vert\,}

\def\ud{\mathrm{d}}
\def\um{\mathrm{m}}
\def\diff {\operatorname{Diff}}
\def\dim{\operatorname{dim}}

   \def\TT{{\mathbb T}}

\def\La{\Lambda}

\def\cA{\mathcal{A}}  \def\cG{\mathcal{G}} \def\cM{\mathcal{M}} 
\def\cB{\mathcal{B}}  \def\cH{\mathcal{H}}  \def\cT{\mathcal{T}}
\def\cC{\mathcal{C}}   \def\cO{\mathcal{O}} \def\cU{\mathcal{U}}
   \def\cP{\mathcal{P}} \def\cV{\mathcal{V}}
    \def\cW{\mathcal{W}}
\def\cF{\mathcal{F}}   \def\cR{\mathcal{R}}

\begin{document}
\vspace{-2cm}
\title{Periodic measures and partially hyperbolic homoclinic classes}
\author{Christian BONATTI and Jinhua ZHANG}
\vspace{-2cm}
\maketitle

\begin{abstract}
 In this paper, we give a precise meaning to the following fact, and we prove it: $C^1$-open and densely, all
the non-hyperbolic ergodic measures generated by a robust cycle  are
approximated by periodic measures.

We apply our technique to  the global setting of partially hyperbolic diffeomorphisms with one dimensional center.
When both strong stable and unstable foliations are minimal,  we get that the closure of the set of  ergodic measures is the union of two
convex sets corresponding to the two possible $s$-indices;   these two convex sets intersect along the closure of the set
 of  non-hyperbolic ergodic measures.
That is the case for robustly transitive perturbation of the time one map of a transitive Anosov flow, or of the skew product
of an Anosov torus diffeomorphism by rotation of the circle.
\end{abstract}

\section{\textbf{Introduction}}\label{sec:1}
\subsection{General setting}

The dynamics of the (uniformly) hyperbolic basic set is considered as well understood, from the topological point of view
as well as from the stochastic point of view. This comes in particular from the existence of Markov partitions which allows us to code  the orbits by
itineraries on which the dynamics acts as a subshift of finite type.  Moreover, these itineraries can be chosen arbitrarily, allowing
to build orbits and measures with a prescribed behavior. One classical result  (Sigmund, 1970) is that any invariant probability measure
$\mu$ supported on a hyperbolic basic set is
accumulated in the weak$*$-topology, by \emph{periodic measures} (i.e. atomic measures associated to a periodic orbit) whose support tends to the
support of $\mu$ for the Hausdorff distance. In particular, the closure of the periodic measures is the  whole convex set
(more precisely the Poulsen simplex) of the invariant probability measures.
This property extends to dynamics with the \emph{specification properties}.

The dynamics on a non-trivial hyperbolic basic set is chaotic.
Sigmund's theorem tells us that any complicated behavior can be described by the simple one in the sense of measure.
For example, we know that any  $C^{2}$ volume preserving Anosov diffeomorphism is ergodic and the whole manifold
is a hyperbolic basic set;  Sigmund's theorem implies that this volume  (whose support is the whole manifold) can
be approximated by periodic measures which are only supported on finite points.

On the other respect, S. Smale and R. Abraham \cite{AS} show that hyperbolic diffeomorphisms
are not dense among $C^1$ diffeomorphisms on a 4-dimensional manifold, and C. Simon~\cite{Si} gives the first such kind of
 example on a 3-manifold.
From the topological point of view, \cite{Sh} and \cite{M1} gave examples of robustly transitive non-hyperbolic diffeomorphisms
 on $\mathbb{T}^4$ and on $\mathbb{T}^{3}$ respectively, then \cite{BD1} gave more general examples of robustly transitive
 non-hyperbolic diffeomorphisms.

The aim of this paper is to investigate to  what extent these properties can be extended in a non-hyperbolic setting.
Many results have already been obtained:

In 1980, A. Katok \cite{Ka} proves that any hyperbolic ergodic measure is  accumulated by periodic measures. 
This result requires  the $C^{1+\alpha}$ setting, and is wrong
in the general 
 $C^1$ setting (see \cite{BCS}). 
 Nevertheless, this result has been extended by \cite{C} in the $C^1$-setting in the case  where the stable/unstable
splitting of a hyperbolic  ergodic measure is a \emph{dominated splitting}.
Even in the smooth setting it is not true in general that the closure of the ergodic measures supported on a homoclinic class is convex
(see for instance counter examples  in  \cite{BG}).
\begin{Question*} Consider a $C^r$-generic diffeomorphism $f$ and a (hyperbolic) periodic point $p$ of $f$. Is it true that the closure of the set of
ergodic measures supported on the homoclinic class $H(p,f)$ is convex?
\end{Question*}
This question remains open in any $C^r$-topology, $r\geq 1$. Let us mention that,  in the $C^1$-generic setting,
 \cite[Theorem 3.10 and Proposition 4.8]{ABC} show  that the closure of the set of periodic measures supported
 on a homoclinic class of a $C^1$-generic diffeomorphism is convex.

 From the stochastic point of view, it was  natural to ask: is it possible to build diffeomorphisms which
 robustly have non-hyperbolic ergodic measures?  \cite{KN} gave the first example of  the robust existence of non-hyperbolic
 ergodic measures, for some partially hyperbolic
 diffeomorphism  whose center foliation is a circle bundle. Recently \cite{BBD2} proved that the existence of non-hyperbolic measures is indeed
 $C^1$-open and dense in the set of diffeomorphisms having a robust cycle.

 Our main results consist in a precise analysis of the periodic orbits whose existence is implied by a robust cycle as the ones considered in \cite{BBD2}.
Using a shadowing lemma in \cite{G1}, we show that every non-hyperbolic measure in a neighborhood of the robust cycle is accumulated by hyperbolic periodic measures.  Putting   the  criteria of \cite{GIKN}
 together with the shadowing lemma in ~\cite{G1},  we  also show  that  the closure
 of the set of these  non-hyperbolic ergodic measures is convex.

\section{Precise statement of the results}

\subsection{Partial hyperbolicity and  robust transitivity}
Our main technical results consist in a local analysis of the periodic orbits appearing in a robust cycle. We present here the consequence of this local analysis
in the setting of robustly transitive partially hyperbolic diffeomorphisms with one dimensional center. In order to present our results,  we need to define carefully this setting
and to recall some known results.

We say that a diffeomorphism $f$  of a compact manifold $M$ is \emph{robustly transitive} if there is a $C^1$-open neighborhood $\cU_f$ of $f$ in $\diff^1(M)$ so that any
$g\in \cU_f$ is \emph{transitive}, that is, $g$ admits a point whose orbit is dense in $M$.  We denote by $\cT(M)\subset \diff^1(M)$, the set of robustly transitive diffeomorphisms. It has been shown by M. Shub~\cite{Sh}  and R.  Ma\~n\'e~\cite{M1} that the  sets of robustly  transitive non-hyperbolic diffeomorphisms  are not empty on $\T^4$ and on $\T^3$ respectively. By definition,
$\cT(M)$ is a $C^1$-open subset of $\diff^1(M)$.

Let $f\colon M\to M$ be a diffeomorphism on a compact manifold $M$, and let $K\subset M$ be a compact invariant subset.  In the literature, there are
several different notions of partial hyperbolicity.  In this paper,  we   say that $K$ is \emph{partially hyperbolic with $1$-dimensional center} if the tangent bundle
$TM|_K$ of $M$ restricted to $K$ admits a $Df$-invariant splitting
$$T_xM= E^{ss}(x)\oplus E^c(x)\oplus E^{uu}(x),$$
where:
\begin{itemize}
 \item the dimensions of $E^{ss}(x)$, $E^{uu}(x)$ are strictly positive and independent of $x\in K$ and $\dim (E^c(x))=1$;
\item there is a Riemannian metric $\|.\|$ on $M$ so that, for any $x\in K$ and for any unit vectors $u\in E^{ss}(x)$, $v\in E^c(x)$ and $w\in E^{uu}(x)$ one has:
$$\|Df(u)\|<\inf\{1,\|Df(v)\|\}\leq \sup\{1,\|Df(v)\|\}<\|Df(w)\|.$$
\end{itemize}
In many usual definitions,  the inequality of the second item only holds for $f^N$ for $N>0$ large enough.  The fact that one can choose the metric so that $N=1$ is due to \cite{Go}.
When $K =M$, one says that $f$ is \emph{partially hyperbolic with one dimensional center bundle}.


 We denote by $\cP\cH(M)\subset \diff^1(M)$ the set of partially hyperbolic diffeomorphisms
with one dimensional center bundle. The partially hyperbolic structure is robust, hence $\cP\cH(M)$ is an open subset of $\diff^1(M)$.

Let $\cO(M)\subset \diff^1(M)$ denote the set of partially hyperbolic (with one dimensional center), robustly transitive diffeomorphisms on $M$. In other words
$$\cO(M)=\cT(M)\cap \cP\cH(M).$$
It has been shown in ~\cite{BV} that there exists a robustly transitive diffeomorphism which has no uniformly hyperbolic invariant subbundles, hence, in general, the set $\cT(M)$ is not contained in $\cP\cH(M)$.

\subsection{Strong stable and strong unstable foliations }
Recall that for any partially hyperbolic diffeomorphism,  there are uniquely defined invariant foliations $\cF^{ss}$ and $\cF^{uu}$ tangent to $E^{ss}$ and to $E^{uu}$ respectively. Our results will depend strongly
on the topological properties of these foliations, and more specifically on the minimality of these foliations. We recall now known results on this aspect.

Recall that
a foliation $\cF$ is \emph{minimal} if every leaf of $\cF$ is  dense.  We say that the strong stable foliation of $f$ is \emph{robustly minimal} if the strong stable foliation
of any diffeomorphism $g$ $C^1$-close enough to $f$ is minimal.

\begin{Remark}\label{r.transitive and mixing}
If the strong stable foliation $\cF^{ss}$ of a diffeomorphism $f\in \cP\cH(M)$ is minimal then $f$ is transitive $($and indeed topologically mixing$)$.
\end{Remark}

We denote by $\cU^s(M)$
(resp. $\cU^u(M)$)  the subset  of
$\cP\cH(M)$ for which the strong stable (resp. strong unstable) foliation is robustly minimal.  We denote $\cU(M)=\cU^s(M)\cap\cU^u(M)$ the set of diffeomorphisms in $\cP\cH(M)$ whose both
 strong stable and unstable foliations are robustly  minimal.  By definitions,  $\cU^s(M)$, $\cU^u(M)$ and $\cU(M)$ are $C^1$-open subsets of $\cP\cH(M)$ and therefore of $\diff^1(M)$.

According to Remark~\ref{r.transitive and mixing}, the sets  $\cU^s(M)$ and $\cU^u(M)$ are contained in $\cT(M)$ and therefore in $\cO(M)$.

\begin{theorem}[\cite{BDU,HHU}] The union  $(\cU^s(M)\cup \cU^u(M))$ is $C^1$-dense (and open)  in $\cO(M)$.
\end{theorem}
In other words, open and densely in the
 setting of partially hyperbolic with $1$-dimensional center and
 robustly transitive diffeomorphisms, one of the strong foliations  is robustly minimal.



\subsection{Center Lyapunov exponent and index of hyperbolic measures}
Consider a  diffeomorphism $f\in\cP\cH(M)$,  denote by  $i=\dim(E^{ss})$,  and recall that $\dim(E^c)=1$. We denote
by $\mathcal{M}_{inv}(f)$ and $\mathcal{M}_{erg}(f)$ the sets of $f$-invariant measures and of  $f$-ergodic measures respectively.

Let $\mu\in\mathcal{M}_{erg}(f)$, then the center Lyapunov exponent of $\mu$ is
$$\lambda^c(\mu)=\int \log \|Df|_{E^c}\| d\mu.$$
One says that the ergodic measure $\mu$ is \emph{hyperbolic} if $\lambda^c(\mu)\neq 0$ and
 $\mu$ is called \emph{non-hyperbolic} if $\lambda^c(\mu)=0$.
The \emph{$s$-index} of a hyperbolic ergodic measure $\mu$
is the number of negative Lyapunov exponents of $\mu$.
In our setting, the $s$-index of a  hyperbolic ergodic  measure  $\mu$ is
\begin{displaymath}
ind(\mu)=\left\{\begin{array}{ll}
i &\textrm{if $\lambda^c(\mu)>0$}
\\
i+1 &\textrm{if $\lambda^c(\mu)<0$ }
\end{array}\right..
\end{displaymath}

We extend by the same formula of the center Lyapunov exponent to the set of all (even non ergodic) invariant probability measures, that is
$$\lambda^c\colon \mathcal{M}_{inv}(f)\to\mathbb{R}, \quad \lambda^c(\mu)=\int \log\| Df|_{E^c}\| d\mu,$$
and we call it the \emph{mean center Lyapunov exponent of $\mu$}. This gives a continuous function in the weak$*$-topology on the space $\mathcal{M}_{inv}(f)$.

  We  denote by  $\mathcal{M}_i(f)$, $\mathcal{M}^{*}(f)$  and $\mathcal{M}_{i+1}(f)$  the sets of
  ergodic measures $\mu$ with positive, vanishing and negative center Lyapunov exponents, respectively.

  We denote by $\mathcal{M}_{Per}(f)$   the subset of $\mathcal{M}_{inv}(f)$ consisting in measures supported on $1$ periodic orbit of $f$.
  We denote by $\mathcal{M}_{Per, i}(f)$ and  $\mathcal{M}_{Per, i+1}(f)$ the subsets of  $\cM_{Per}(f)$, consisting in  measures
  supported on a single hyperbolic periodic orbit   of index $i$ and $i+1$ respectively.

  Finally, if $\mathcal{E}\subset \mathcal{M}_{inv}(f)$, then $\overline{\mathcal{E}}$ denotes its weak$*$-closure.

\subsection{Main results  when both strong foliations are minimal}

Let $M$ be a compact manifold.
We denote by  $\mathcal{V}(M)\subset \diff^1(M)$  the set of  diffeomorphisms such that:
 \begin{itemize}\item $f\in \cU(M)$ , that is, $f$ is partially hyperbolic with one dimensional center bundle and has both strong stable and unstable foliations which are robustly minimal.
 \item $f$ admits hyperbolic periodic points $p_f$ of $s$-index $i$ and $q_f$ of $s$-index $i+1$.
 \end{itemize}

 \begin{Remark}
 $\mathcal{V}(M)$ is, by definition,  an open subset of $\diff^1(M)$. \cite{BDU} shows that there are compact $3$-manifold $M$ for which $\cV(M)$ is not empty:
 the time one map of a transitive Anosov flow on a manifold $M$ admits a  smooth perturbation in $\mathcal{V}(M)$; the same occurs in  the skew product of
   linear Anosov automorphisms of the torus $\TT^2$ by rotations of the circle.
\end{Remark}

For any $f\in \cV(M)$ and any hyperbolic periodic point $x$, one has the following properties:
\begin{itemize}

\item the minimality of both foliations implies that the manifold $M$ is the whole homoclinic class of $x$. Furthermore any
two hyperbolic periodic points of   same index are homoclinically related.  This implies that the closure of the set of  the hyperbolic periodic measures of a given index is convex.  In other words
$\overline{\cM_{Per,i}(f)}$ and $\overline{\cM_{Per,i+1}(f)}$ are convex.

 \item clearly $\cM_{Per,i}(f)\subset \cM_i(f)$ and $\cM_{Per,i+1}(f)\subset \cM_{i+1}(f)$.
 According to \cite{C}, in this partially hyperbolic setting,  every hyperbolic ergodic measure is weak$*$-limit of
 periodic measures of the same index.  One gets therefore:
 $$\overline{\cM_{Per,i}(f)}=\overline{\cM_i(f)} \mbox{ and } \overline{\cM_{Per,i+1}(f)}=\overline{\cM_{i+1}(f)}$$
\end{itemize}

\begin{Remark} \label{r.generic}
 As a direct consequence of \cite{M2,ABC}, there is a $C^1$-residual subset $\cG$ of $\cV(M)$ so that for every $f\in\cG$,
   every invariant (not necessarily ergodic) measure is the weak$*$-limit of periodic measures.
 In particular,  $\cM_{inv}(f)$ is a Poulsen simplex. 
More precisely, \cite{M2} implies that generically any ergodic measure is the limit of periodic measures,
and \cite{ABC} shows that, for generic robustly transitive diffeomorphisms the closure of the periodic measures is convex.
\end{Remark}

We don't know if the generic properties stated in Remark~\ref{r.generic} hold for every $f\in\cV(M)$. The aim of our main result in this setting is to recover most
of these properties for a $C^1$ open and dense subset
of $\cV(M)$.
Recall that $\cM^{*}(f)$,  $\cM_i(f)$ and $\cM_{i+1}(f)$ are  the sets of ergodic measures with vanishing, positive and negative center Lyapunov exponents respectively.

  \begin{mainthm} \label{thmA} Let  $M$ be a closed manifold.
  There exists a  $C^1$ open and dense subset $\widetilde{\mathcal{V}}(M)$ of $\mathcal{V}(M)$ such that for any $f\in\widetilde{\mathcal{V}}(M)$, we have the followings:
  \begin{enumerate}
  \item $$\begin{array}{ccl}
           \overline{\mathcal{M}^*(f)}&=&\overline{\mathcal{M}_{i}(f)}\cap\{\mu\in \mathcal{M}_{inv}(f),\lambda^c(\mu)=0\}\\
           &=&\overline{\mathcal{M}_{i+1}(f)}\cap\{\mu\in \mathcal{M}_{inv}(f), \lambda^c(\mu)=0\}\\
 &=&\overline{\mathcal{M}_{i}(f)}\cap\overline{\mathcal{M}_{i+1}(f)}.
  \end{array}$$
  In particular, the closure $\overline{\mathcal{M}^*(f)}$ of the non-hyperbolic ergodic measures is convex. Furthermore, every non-hyperbolic
  ergodic measure is the weak$*$-limit of hyperbolic periodic orbits of both indices $i$ and $i+1$.

\item There exist two  compact $f$-invariant $($uniformly$)$ hyperbolic   sets $K_i\subset M$ and $K_{i+1}\subset M$ of $s$-index $i$ and $i+1$  respectively, with the following property:
  for any $\mu\in\overline{\mathcal{M}_i(f)}$ (resp. $\overline{\mathcal{M}_{i+1}(f)}$), there exists an invariant measure $\nu$  supported on $K_{i+1}$ (resp. $K_i$) such that
  the segment $\{\alpha\mu+(1-\alpha)\nu|\,\alpha\in[0,1]\}$ is contained in the closure of the set of periodic measures.
  \end{enumerate}
  \end{mainthm}
  \begin{Remark} The existence of non-hyperbolic ergodic measure is guaranteed by \cite{BBD2} (see also ~\cite{BZ}).
\end{Remark}

The proof of this theorem is based on the semi-local setting (the next subsection).

Theorem \ref{thmA} shows that the closure of the set of ergodic measures for $f\in\tilde{\cV}(M)$
is the union of two convex sets $\overline{\cM_{Per,i}(f)}$ and $\overline{\cM_{Per,i+1}(f)}$, which intersect along $\overline{\cM^{*}(f)}$. The last item of the theorem shows that the union of these convex sets ``is not far from being convex'',
but we did not get the convexity. In other words, we don't know if any $f$-invariant measure is accumulated by ergodic measures.
\begin{Question}Does there exist an open dense subset of $\cV(M)$ such that periodic measures are dense among invariant measures?
\end{Question}

\subsubsection{The $C^1$-generic case}
As told before in Remark~\ref{r.generic},  for $C^1$-generic $f$ in the set $\tilde{\cV}(M)$, the closure  $\overline{\cM_{Per}(f)}$ of the set of hyperbolic periodic measures is convex, and coincides with
 the set $\cM_{inv}(f)$ of all invariant measures.

Our result implies:
\begin{corollary} For $C^1$-generic  $f\in\tilde{\cV}(M)$, every invariant (a priori non ergodic) measure whose
mean center Lyapunov exponent vanishes is approached by non-hyperbolic ergodic measures.
In formula:
$$ \{\mu\in\cM_{inv}(f), \lambda^c(\mu)=0\}= \overline{\cM^*(f)} .$$

\end{corollary}

It is well known that the decomposition of  an invariant  measure in the convex sum of ergodic measures is unique and one calls  $\cM_{inv}(f)$  a \emph{ Choquet simplex}. As we split $\cM(f)$  into  several convex sets (the ones with
positive, vanishing and negative center Lyapunov exponents respectively) it is natural to ask if these sets are Choquet simplices too. We could not answer to this question
in the whole general situation,  but in the $C^1$-generic setting we can give a negative answer:

\begin{proposition}\label{p.non simplex}For $C^1$-generic $f$ in $\tilde{\cV}(M)$,
none of the three compact convex sets $\overline{\cM_i(f)}$,
$\overline{\cM_{i+1}(f)}$ and $\overline{\cM^*(f)}$ is a Choquet simplex.
\end{proposition}

\subsection{Main results  when only one strong foliation is minimal}
As we mentioned before, it is known  that open and densely in the set $\cO(M)$ of robustly transitive, partially hyperbolic (with one dimensional center) diffeomorphisms, one of the strong   foliations is robustly minimal: $\cU^{s}(M)\cup\cU^{u}(M)$ is dense in $\cO(M)$. As far as we know there are no known examples where $\cU(M)=\cU^{s}(M)\cap\cU^{u}(M)$ is not dense in $\cO(M)$. Nevertheless,  we cannot discard this possibility.

In general, when only one of the strong foliations is minimal, we don't know if the non-hyperbolic ergodic measures are accumulated by periodic orbits, and if the closure of the set of non-hyperbolic ergodic measures is convex. Nevertheless,  there is an important example where we could recover these properties.

   R. Ma\~n\'e \cite{M1} gave an example that for linear Anosov diffeomorphisms on $\mathbb{T}^3$ of $s$-index $1$ with three ways dominated splitting, one can do $DA$ to get an open subset $\mathcal{W}$ of $\diff^1(\mathbb{T}^3)$ where all the diffeomorphisms are non-hyperbolic and transitive (see the precise definition of $\cW$ in Section~\ref{s.mane example}).
     \cite{BDU} proved that   the robustly transitive diffeomorphism in $\cW$ has minimal strong stable foliation (see also \cite{PS}), that is,
   $\cW$ is contained in $\cU^{s}(\mathbb{T}^{3})$.

\begin{mainthm}\label{thmB} There exists an open and   dense subset $\tilde{\mathcal{W}}$ of $\mathcal{W}$ such that for any $f\in\tilde{\mathcal{W}}$, one has
 \begin{itemize}
 \item Any non-hyperbolic ergodic measure is approximated by periodic measures of s-index $1$ (recall that all the periodic orbits of s-index $1$ are homoclinically related);
 \item The closure  $\overline{\cM^{*}(f)}$ of the set of non-hyperbolic ergodic measures is convex.
\end{itemize}
\end{mainthm}

Since the strong stable foliation of every $f\in \cW$  is minimal,  all hyperbolic periodic orbits of $s$-index $1$ are
homoclinically related. By  transitivity, one can show that the unstable manifold of hyperbolic periodic orbit of $s$-index $1$ is dense on the manifold. Hence  $M$ is the homoclinic class of every periodic orbit of $s$-index $1$.  By~\cite[Theorem E]{BDPR}, for an open and dense subset of $\cW$, the manifold $M$ is also the homoclinic class of a periodic orbit of $s$-index $2$.
 \begin{Question} For Ma\~n\'e's example, given two hyperbolic periodic orbits $Q_1$ and $Q_2$ of s-index $2$, are $Q_1$ and $Q_2$ homoclinically related?
 \end{Question}
 \begin{Remark}    If the answer to Question 2 is yes, although we can not get the whole convexity of the set of periodic measures,  one can show that for Ma\~ n\'e 's  example, the set of  ergodic measures  is path connected $($for the definition  see~\cite{Sig2}$)$.
 \end{Remark}

 \subsection{Some results in the semi-local setting of  robust cycles}
We start with our assumption in the semi-local setting without technical definitions and the precise definition of the terminology we use here would be given in   the next section.

Consider $f\in\diff^1(M)$. Let  $(\Lambda,U,\mathcal{C}^{uu},\mathfrak{D})$ be a  \emph{blender horseshoe} of u-index $i+1$  and
 $\mathcal{O}_{q}$ be a hyperbolic periodic orbit of u-index $i$.
 We assume that $\Lambda$ and $\cO_q$ form  a special  robust cycle called  \emph{split-flip-flop configuration}.
We fix a  small neighborhood $V$ of the split flip-flop configuration so that the maximal invariant set $\tilde \La$ in the closure $\bar V$ admits a partially hyperbolic splitting $E^{ss}\oplus E^c\oplus E^{uu}$ with $\dim (E^c)=1$. Assume, in addition, that there exists a $Df$-strictly invariant center unstable cone field $\cC^u_V$ which is a continuous extension of the center unstable cone field $\cC^u$ in $U$.

In ~\cite{ABC}, it has been shown that in the $C^1$ generic setting, for each homoclinic  class, the closure of the set of periodic measures is convex. Hence, for each  partially hyperbolic homoclinic  class with center dimension one, under $C^1$-generic setting, the closure of the  set of hyperbolic ergodic measures
 is convex.
  Here, in the $C^1$ open setting, we get some kind of `convexity' and it would be used for the global setting.
  \begin{mainthm}\label{thmE}Under the assumption above. There exists an invariant measure $\mu\in\mathcal{M}_{inv}(\Lambda, f)$ such that the segment $\{\alpha\mu+(1-\alpha)\delta_{\mathcal{O}_q}, \alpha\in[0,1]\}$ is contained
in the closure of the set of periodic measures whose support are inside $V$.
\end{mainthm}

   In any small neighborhood of the split flip flop configuration,  the existence of non-hyperbolic ergodic measure approached
   by periodic measures has been proven in~\cite{BZ}.
    Conversely, we can prove that 
    the non-hyperbolic ergodic  measures supported in a small neighborhood of the robust cycle  are  approximated by hyperbolic periodic measures. To be precise:
   \begin{mainthm}\label{thmF} Under the assumption above. There exists a small neighborhood $V_0\subset V$ of the split flip flop configuration such that  for any non-hyperbolic ergodic  measure $\nu$ supported on the maximal invariant set $\tilde{\La}_0$ in $V_0$, there exists a sequence of periodic orbits $\{\cO_{p_{n}}\}_{n\in\N}$ which are homoclinically related to $\Lambda$ such that $\delta_{\mathcal{O}_{p_{n}}}$ converges to $\nu$.
   \end{mainthm}
   \begin{Remark}
    \begin{enumerate}
    \item If the support of $\nu$ intersects the boundary of $V_0$, the sequence of periodic orbits we find might  intersect the complement of  $\tilde{\Lambda}_0$;
        \item the choice of $V_0$ is uniform for the diffeomorphisms in a small  $C^1$  neighborhood of $f$.
    \end{enumerate}
   \end{Remark}

\section{Preliminaries}\label{s.Preliminaries}

  In this section, we will collect some notations and some results that we need.  We start by recalling very classical notions, as Choquet simplex, Poulsen simplex and dominated splitting.
  Then we recall our main  tools.  More precisely, our results consist in applying  two tools in a very specific setting.  The tools are:
  \begin{itemize}
  \item the \emph{\cite{GIKN} criterion} for ensuring that the limit measure of periodic measures  is ergodic.

  \item a shadowing Lemma due to Liao \cite{Liao1} and Gan   \cite{G1} : this will allow us to prove the existence periodic orbits with a prescribed itinerary.
   An important consequence presented in \cite{C} shows that any hyperbolic ergodic measure is  approached by  periodic orbits, under a dominated setting.
  \end{itemize}
  Our setting will be  a specific robust cycle   called \emph{split flip flop configuration}. The main interest of the split  flip flop configuration is
  that it  appears  open and densely in the setting of robust cycle.

\subsection{Dominated splitting and hyperbolicity}
 Recall that a $Df$-invariant splitting $T_{K}M=E\oplus F$ over a compact $f$-invariant set $K$ is a \emph{dominated splitting}, if there exists
 $\lambda\in(0,1)$ and a metric $\|\cdot\|$ such that
 $$ \norm{Df |_{E(x)}}\cdot \norm{Df ^{-1}|_{F(f(x))}}< \lambda,\textrm{ for any $x\in K$.}$$

A compact $f$-invariant set $K$ is called  a \emph{   hyperbolic set},  if there exists an invariant hyperbolic splitting  $T_{K}M=E^s\oplus E^{u}$, that is,   $E^s$ is uniformly contracting and $E^{u}$ is uniformly expanding under $Df$.
A hyperbolic set $K$ is called a \emph{hyperbolic basic set of s-index $i$} if one has
\begin{itemize}
\item $K$ is transitive and $\dim (E^s)=i$;
 \item there exists an open neighborhood $U$ of $K$ such that $K$ is the \emph{maximal invariant set in $U$}, that is,
$$K=\cap_{i\in\Z}f^{i}(U).$$
\end{itemize}
We denote by $ind(\La)$ the s-index of $\La$.

 One    important property of hyperbolic basic sets in the sense of measure is the following theorem:
 \begin{theorem}\label{sigmund}\cite[Theorem 1]{Sig1}  Let $f\in\diff^1(M)$ and $\Lambda$ be a hyperbolic basic set. Then any $f$ invariant measure supported on $\La$ is approximated by periodic measures.
 \end{theorem}
\subsection{Homoclinic class }
\begin{definition} Let $f\in\diff^1(M)$. Given two hyperbolic periodic orbits   $\mathcal{O}_p$ and $\mathcal{O}_q$ of $f$.
    $\mathcal{O}_p$ and $\mathcal{O}_q$ are said to be  \emph{homoclinically related}, if  there exist a non-empty transverse intersection between $W^s(\mathcal{O}_p)$ and $W^u(\mathcal{O}_q)$, and a non-empty transverse intersection between $W^u(\mathcal{O}_p)$ and  $W^s(\mathcal{O}_q)$.
\end{definition}
Let $\mathcal{O}_p$ be a hyperbolic periodic orbit, the \emph{homoclinic class of $\mathcal{O}_p$} is defined as:
$$H(p,f):=\overline{\{\mathcal{O}_q|\textrm{  The hyperbolic perodic orbit $\mathcal{O}_q$ is homoclinically  related to $\mathcal{O}_p$}\}}.$$

Let  $\mathcal{O}_p$ and $\mathcal{O}_q$ be  two hyperbolic periodic orbits and  $V$ be an open  neighborhood of $\mathcal{O}_p\cup\mathcal{O}_q$.
 We say that $\mathcal{O}_p$ and $\mathcal{O}_q$ are \emph{homoclinically related inside $V$},
 if there exist two transverse intersections  $x\in W^{u}(\cO_p)\cap W^{s}(\cO_q)$ and  $y\in W^{s}(\cO_p)\cap W^{u}(\cO_q)$
 such that $\overline{Orb(x)}\cup \overline{Orb(y)}\subset V.$
\subsection{Lyapunov exponents, Oseledets splitting and  hyperbolic ergodic measure}
In the celebrated paper \cite{O}, V. Oseledets proves that for any ergodic measure $\mu$ of a diffeomorphism $f$, we have the following:
\begin{enumerate}\item there exists a $\mu$-full measure set $K$ such that $f(K)=K$;
\item there exist $s\leq \dim(M)$ numbers $\lambda_1<\cdots<\lambda_s$ and an invariant measurable splitting over $K$ of the form $T_{K}M=E_1\oplus\cdots\oplus E_s$ such that for any integer $1\leq t\leq s$, any $x\in K$ and any $v\in \oplus_{i=1}^{t}E_i(x)\backslash \oplus_{i=1}^{t-1}E_i(x)$, we have that
     $$\lim_{n\rightarrow+\infty}\frac{1}{n}\log{\norm{Df^n_{x}v}}=\lambda_t.$$
\end{enumerate}
The numbers $\lambda_1,\cdots,\lambda_s$ are called the \emph{Lyapunov exponents} of $\mu$, the full measure set $K$ is called the \emph{Oseledets basin} of $\mu$ and the splitting $T_{K}M=E_1\oplus\cdots\oplus E_s$ is called the \emph{Oseledets splitting} of $\mu$.

An ergodic measure $\mu$ is called a \emph{hyperbolic ergodic measure}, if all the Lyapunov exponents of $\mu$ are non-zero. Let $K$ be the Oseledets basin of $\mu$, we denote by $E^-\oplus E^{+}$   the invariant splitting over $K$ such that all the Lyapunov exponents along $E^-$ are negative and all the Lyapunov exponents along $E^{+}$ are positive. Then the invariant splitting $T_KM=E^-\oplus E^{+}$
is called the \emph{non-uniform hyperbolic splitting}.
We say that the \emph{non-uniform hyperbolic splitting is dominated} if there exists a dominated splitting over the closure of $K$ of the form $T_{\overline{K}}M=E\oplus F$ such that $\dim(E)=\dim(E^-)$.

 \subsection{Choquet and Poulsen Simplex}

  \begin{definition} Let $K$ be a non-empty compact   convex subset of a locally convex vector space. Then
  \begin{itemize}
  \item $K$ is said to be a \emph{Choquet simplex}, if   every point of $K$ is the barycenter of a unique probability measure supported on the set of extreme points of $K$.
  \item $K$ is said to be  a \emph{Poulsen simplex} if $K$ is a Choquet Simplex  so that  the set of extreme points of $K$ is strictly contained in $K$ and is dense in $K$.
  \end{itemize}

  \end{definition}
Let $C$ be a compact  $f$-invariant set. We denote by $\mathcal{M}_{inv}(C,f)$ the set of invariant measures supported on $C$.
 A classical result  is that $\mathcal{M}_{inv}(C,f)$ is always a Choquet simplex, hence
  Sigmund's theorem (see \cite{Sig1}) shows  that if $C$ is  a hyperbolic basic set then $\cM_{inv}(C,f)$ is a  Poulsen Simplex.
  \subsection{Good approximation and \cite{GIKN}  criterion}
  In this subsection, we  state the ~\cite{GIKN} criterion ensuring that a sequence of periodic measures converges to an ergodic measure.
  This criterion is firstly used in \cite{GIKN,KN} and is  developed in \cite{BDG} for building  non-hyperbolic ergodic measures as the limit of periodic measures.
  \begin{definition}\label{good for}Given a compact metric space  $(X,\ud)$  and let  $f:X\mapsto X$ be a continuous map.  Fix $\epsilon>0$ and $\kappa\in(0,1)$.  Let $\gamma_1$ and $\gamma_2$ be two periodic orbits of $f$. Then,
  the periodic orbit $\gamma_{1}$ is said to be a  \emph{$(\epsilon,\kappa)$ good  approximation} for   $\gamma_{2}$,
   if there exist  a subset $\gamma_{1,\epsilon}$ of $\gamma_{1}$ and a projection $\cP :\gamma_{1,\epsilon}\rightarrow\gamma_{2}$  such that:
  \begin{itemize}
  \item  for every $y\in\gamma_{1,\epsilon}$ and every $i=0,\cdots,\pi(\gamma_{2})-1$, one has
   $$\ud(f^{i}(y),f^{i}(\cP(y)))<\epsilon;$$
  \item the proportion of $\gamma_{1,\epsilon}$ in $\gamma_1$ is larger that $\kappa$.  In formula:
      $$\frac{\#\gamma_{1,\epsilon}}{\pi(\gamma_{1})}\geq\kappa.$$
  \item  the cardinal of the pre-image $\cP^{-1}(x)$  is the same for all $x\in\gamma_{2}$.
  \end{itemize}
  \end{definition}
  We can now state the \cite{GIKN} criterion refined in \cite{BDG}:
  \begin{lemma}\label{limit} \cite[Lemma 2.5]{BDG} Let $(X,d)$ be a compact metric space and $f:X\mapsto X$ be a homeomorphism. Let $\{\gamma_{n}\}_{n\in\N}$ be a sequence of periodic orbits
  whose periods tend to infinity.  We denote by  $\mu_{n}$   the Dirac measure of $\gamma_{n}$.

  Assume that the orbit $\gamma_{n+1}$ is a $(\epsilon_{n}, \kappa_{n})$ good approximation  for $\gamma_{n}$, where  $\epsilon_{n}>0$ and
   $0<\kappa_{n}<1$ satisfy
  \begin{displaymath} \sum_{n}\epsilon_{n}<\infty \textrm{ and } \prod_{n} \kappa_{n}>0.
  \end{displaymath}

  Then the sequence $\mu_{n}$ converges to an ergodic measure $\nu$ whose support is given by
  \begin{displaymath} supp\,\nu=\cap_{n=1}^{\infty}\overline{\cup_{k=n}^{\infty}\gamma_{k}}.
  \end{displaymath}
\end{lemma}

\subsection{Shadowing lemma and approaching hyperbolic measures by periodic orbits}
In this paper we don't construct any periodic orbits by perturbations.
We find these periodic orbits by  a shadowing lemma which is firstly given by S. Liao \cite{Liao1} and is developed by S. Gan~\cite{G1}.

    Let $\Lambda$ be a compact  $f$-invariant  set,  admitting a dominated splitting of the form  $T_{\Lambda}M=E\oplus F$.
      For any $\lambda\in(0,1)$, an orbit segment $\{x,n\}:=\{x,\cdots,f^{n}(x)\}$ contained in $\Lambda$ is called a \emph{$\lambda$-quasi hyperbolic string}, if the followings are satisfied:
      \begin{itemize}
       \item {\bf Uniform contraction of $E$ by $Df$, from $x$ to $f^n(x)$:}
       $$\prod_{i=0}^{k-1}\|Df|_{E(f^{i}(x))}\|\leq\lambda^{k},$$ for every  $k=1,\cdots,n$;
       \item {\bf Uniform contraction of $F$ by $Df^{-1}$, from $f^n(x)$ to $x$}         $$\prod_{i=k}^{n-1}\um(Df|_{F(f^{i}(x))})\geq(\lambda^{-1})^{n-k},$$ for every $k=0,\cdots,n-1$;
      \end{itemize}

 \begin{Remark}
By definition, one can   check  that a $\lambda$-quasi hyperbolic string is also a $\frac{1+\lambda}{2}$-quasi hyperbolic string.
\end{Remark}




Now, we state the \emph{shadowing lemma for quasi hyperbolic pseudo orbit}:
 \begin{lemma} \label{shadow}\cite{Liao1,G1} Let $f\in\diff^1(M)$ and $\La$ be a compact $f$-invariant set.
   Assume that $\Lambda$  exhibits a  dominated  splitting $T_{\Lambda}M=E\oplus F$.

  Then, for any $\lambda\in(0,1)$, there exist  $L>0$ and $d_{0}>0$ such that for any $d\in(0,d_{0}]$ and any
  $\lambda$-quasi hyperbolic  string $\{x, n\}$ satisfying that $\ud(f^n(x),x)\leq d$,
  there exists a periodic  point $p$ of period $n$ which   shadows $\{x, n\}$ in the distance  of $L\cdot d$,
  that is,  
  $$\ud(f^i(p),f^i(x))<L\cdot d, \textrm{ for any  $i=0,\cdots, n-1$.}$$

  \end{lemma}
We call the orbit segment $\{x,n\}$ above as \emph{the $\lambda$-quasi hyperbolic periodic $d$-pseudo orbit.}

Lemma~\ref{shadow} (together with the $C^1$-Pesin theory  in \cite{ABC})  has been used by S. Crovisier~\cite{C} for approaching hyperbolic ergodic measures with periodic orbits:

\begin{proposition}\label{c1}\cite[Proposition 1.4]{C}
 Let $f\in\diff^1(M)$ and $\mu$ be a hyperbolic ergodic
measure whose non-uniform hyperbolic splitting $E^{-}\oplus E^{+}$ is dominated. Then $\mu$ is supported on a homoclinic class.

Moreover, there exists a sequence of periodic orbits $\{\gamma_n\}_{n\in\N}$ of  $s$-index $\dim(E^{-})$ which are pairwise  homoclinically related,
such that $\gamma_n$ converges to the support of $\mu$ for the Hausdorff topology and  the Dirac measure
supported on $\gamma_n$ converges to $\mu$ in the weak${*}$-topology.
\end{proposition}
\begin{Remark} This  result is firstly obtained by \cite{G2} on surfaces.
\end{Remark}
\subsection{Plaque family   and estimate  on the size of   invariant manifold}
 In this section, let $\La$ be a compact  $f$-invariant set with a dominated splitting $T_{\La}M=E\oplus F$.

  We recall the \emph{Plaque family theorem} by \cite{HPS} showing that there exist invariant plaque families for dominated splitting.
  Given a  continuous bundle $G$ over a   set $K$, for any $x\in K$ and $r>0$, we denote by $G_x(r)=\{v\in G_x|\norm{v}\leq r\}$ and  denote by $G(r)=\cup_{x\in K}G_x(r)$.
  \begin{lemma}\label{l.plaque family} Let $\Lambda$ be a compact $f$-invariant set admitting a dominated splitting $T_{\Lambda}M=E\oplus F$. Then there exist two  continuous maps  $W^{cs}: E(1)\rightarrow M$ and $W^{cu}:  F(1) \rightarrow M$ satisfying the followings:
   \begin{itemize}
   \item for any $x\in\La$, the induced map  $W^{cs}_x: E_x(1)\mapsto M$ (resp. $W^{cu}_x:  F_x(1)\mapsto M$)  is a  $C^1$ embedding and is  tangent to $E_x$ (resp. $F_x$) at the point $x=W^{cs}_x(0)$ (reps. $x=W^{cu}_x(0)$).
   \item the families  $\{W^{cs}_x\}_{x\in\La}$ and $\{W^{cu}_x\}_{x\in\La}$ of $C^1$ embedding maps  are continuous;
   \item there exists a neighborhood $U_{E}$ (resp. $U_F$) of zero section in $E$ (resp. $F$) such that the image of $W^{cs}_x(E_x\cap U_E)$ (resp. $W^{cu}_x(F_x\cap U_F)$) by $f$ (resp. $f^{-1}$) is contained in $W^{cs}_x(E_x)$ (reps. $W^{cu}_x(F_x\cap U_F)$).
   \end{itemize}
   \end{lemma}
   We  denote by  $\mathcal{W}_{\delta}^{cu}(x)=W_x^{cu}(F_x(\delta))$ and $\mathcal{W}_{\delta}^{cs}(x)=W^{cs}_x(E_x(\delta))$, for $\delta\in(0,1]$.

\begin{definition} Given $\lambda\in(0,1)$. A point $x\in\Lambda$ is
 \emph{a $(\lambda, E)$-contracting point,} if we have that
$$\prod_{i=0}^{j-1}\norm{Df|_{E(f^{i}(x))}}\leq \lambda^j, \textrm{ for any positive integer $ j$.}$$
Similarly, we can define the $(\lambda^{-1},F)$-expanding point  which is a $(\lambda, F)$ contracting point  for $f^{-1}$.
\end{definition}
The following lemma guarantees the existence of stable manifolds at the $(\lambda,E)$-contracting points. The proof is classical see for instance \cite[Section8.2]{ABC}. According Lemma \ref{l.plaque family}, we fix the plaque families $W^{cs}$ and $W^{cu}$.
\begin{lemma}\label{lp} For any $\lambda\in(0,1)$, there exists $\eta>0$  such that for any $(\lambda, E)$ contracting point  $x$, we have that  the disc $W^{cs}_{\eta}(x)$ is contained in the stable manifold of $x$.
  \end{lemma}
  \begin{Remark} Similar result holds for  $(\lambda^{-1},F)$-expanding points.
  \end{Remark}
 To find the  $(\lambda,E)$-contracting points, we need the following well known Pliss lemma:
 \begin{lemma}\label{l.Pliss lemma}\cite{P}
  Let $a_{1},\ldots,a_{n}$ be a sequence of numbers bounded from above by a number $b$. Assume that there exists a number
  $c<b$ such that $\sum_{i=1}^{n}a_{i}\geq n\cdot c.$

    Then for any number $c^{\prime}<c$, there exist $l$ integers $i_{1},\ldots,i_{l}\subset [1,n]$ such that
   for any $k=1,\ldots,l$, we have that
    $$\sum_{i=j}^{i_{k}}a_{i}\geq (i_{k}-j+1)c^{\prime}, \textrm{  for any $j=1,\ldots,i_{k}$.}$$
 Moreover, one has that
     $$\frac{l}{n}\geq\frac{c-c^{\prime}}{b-c^{\prime}}.$$
  \end{lemma}

 \subsection{Blender}

   \emph{Blender} is a powerful tool   in the study of robustly non-hyperbolic phenomena.
   In this subsection, we will state  a new definition of blender recently defined by \cite{BBD2},  and a special   blender called \emph{Blender horseshoe} given in \cite{BD2}.

\subsubsection{Dynamically defined blender}
   Let's first  recall some notations in \cite{BBD2}.
  Let $D^{i}(M)$ be  the set of $C^1$
  embedded $i$-dimensional compact
  discs on a  compact Riemannian manifold $M$. We endow $D^{i}(M)$ with $C^1$-topology in the following way: for any $D\in D^{i}(M)$, which is the image of the embedding  $\phi:\mathbb{D}^{i}\mapsto M$
  where $\mathbb{D}^{i}$ is the $i$-dimensional closed unit disc in $\mathbb{R}^{i}$, we define the $C^{1}$ neighborhood of $D$  as the set of the images of all the embedding maps contained in a $C^1$ neighborhood of $\phi$.

  For any $D_1, D_2\in D^{i}(M)$, one can define their distance as
  $$\rho(D_1,D_2)=\ud_{Haus}(TD_1,TD_2)+\ud_{Haus}(T\partial{D_1}, T\partial{D_2}),$$
  where $\ud_{Haus}(\cdot,\cdot)$ denotes the Hausdorff distance on the corresponding Grassmann manifold.
  It has been proven in ~\cite[Section 3.1]{BBD2} that the distance $\rho(\cdot,\cdot)$
  induces the $C^1$-topology in $D^{i}(M)$.

  Let $f\in\diff^{1}(M)$, and $\mathfrak{D}$ be a family of $i$-dimensional embedded discs. For any $\epsilon>0$, we denote by $\mathcal{V}_{\epsilon}(\mathfrak{D})$ the $\epsilon$-neighborhood
  of $\mathfrak{D}$ for the distance $\rho(\cdot,\cdot)$.

 The family $\mathfrak{D}$ is called  a \emph{strictly invariant family}, if there exists $\epsilon>0$ such that for any $D\in\mathcal{V}_{\epsilon}(\mathfrak{D})$, its image $f(D)$ contains an element of $\mathfrak{D}$. The number $\epsilon$
 is called the strength of the strictly invariant family $\mathfrak{D}$.

  \begin{definition}\label{dynamical blender}$(Dynamical\,\, Blender)$ Let $f\in \diff^{1}(M)$. A
  hyperbolic basic set $\Lambda$ of $f$ is called a \emph{dynamically defined cu-blender}
  of uu-index $i$, if the followings are satisfied:
  \begin{itemize}
  \item  there exists a    dominated splitting of the form $T_{\Lambda}M=E^{s}\oplus E^{c}\oplus E^{uu}$ over $\La$;
  where $\dim (E^{s})=ind(\Lambda)$, $\dim (E^{c})>0$ and $\dim (E^{uu})=i$.

      \item there exists a neighborhood $U$ of $\Lambda$ such that $\Lambda=\bigcap_{n\in\mathbb{Z}}f^{n}(U)$ and there exists a
      $Df$-strictly invariant continuous
      cone field $\mathcal{C}^{uu}$ of index $i$ defined on $\overline{U}$;
          \item there is a strictly invariant family $\mathfrak{D}\subset D^{i}(M)$ of discs with strength $\epsilon>0$ such that every disc $D\in        \mathcal{V}_{\epsilon}(\mathfrak{D})$ is tangent to $\mathcal{C}^{uu}$ and is contained in $U$.

            \end{itemize}
The set $U$ is called the  domain of $\Lambda$, the cone field $\mathcal{C}^{uu}$ is called  the strong unstable cone field of $\Lambda$ and
the family $\mathfrak{D}$ is called  strictly invariant family of discs. We denote the dynamically defined cu-blender as $(\La, U,\cC^{uu},\mathfrak{D})$
\end{definition}

  We can also define the \emph{cs-blender} which is a \emph{cu-blender} for the reversed dynamics.
  \begin{Remark}\label{r.domain of blender}
   \cite[Scholium 3.15]{BBD2} Let $(\La, U, \cC^{uu},\mathfrak{D})$ be a dynamically defined cu-blender, there exists a disc in the local strong unstable manifold of a point in $\La$ which are approximated by discs in $\mathfrak{D}$.

  \end{Remark}


The main property of a dynamically defined blender is the following:
    \begin{lemma}\label{robustness of blender} \cite[Lemma 3.14]{BBD2}
     Let $(\Lambda,U,\mathcal{C}^{uu},\mathfrak{D})$ be a dynamically defined cu-blender and  $\epsilon$ be the strength of the strictly invariant family $\mathfrak{D}$.

    Then there exists a $C^1$ neighborhood $\mathcal{U}$ of $f$ such that for any $g\in\cU$, one has
    \begin{itemize}\item  Let  $\Lambda_g$ be the continuation of $\La$. For any $D\in\mathcal{V}_{\epsilon/2}(\mathfrak{D})$, one has that
    $W^{s}(\Lambda_{g})\cap D\neq\emptyset$;
    \item  the continuation  $(\Lambda_g,U,\mathcal{C}^{uu},\mathcal{V}_{\epsilon/2}(\mathfrak{D}))$ is a dynamically defined blender for $g$.
    \end{itemize}
    We call the open family $\mathcal{V}_{\epsilon/2}(\mathfrak{D})$ the superposition region of the blender $(\Lambda,U,\mathcal{C}^{uu},\mathfrak{D})$.
    \end{lemma}
\subsubsection{Blender horseshoe}
In this part, we recall the main feature of a special and simplest  blender called \emph{Blender Horseshoe} (for specific definition see \cite{BD3}).

Consider $\R^{n}=\R^{s}\oplus\R \oplus\R^{u}$.  For $\alpha\in(0,1)$,
  we   define the following cone fields:
 $$\cC_{\alpha}^{s}(x)=\{v=(v^{s},v^{c},v^{u})\in \mathbb{R}^{s}\oplus\mathbb{R}\oplus\mathbb{R}^{u}=T_{x}M| \| v^{c}+v^{u}\|\leq\alpha\|v^{s}\|\}$$
  $$\cC_{\alpha}^{u}(x)=\{v=(v^{s},v^{c},v^{u})\in \mathbb{R}^{s}\oplus\mathbb{R}\oplus\mathbb{R}^{u}=T_{x}M| \| v^{s}\|\leq\alpha\|v^{c}+v^{u}\|\}$$
  $$\cC_{\alpha}^{uu}(x)=\{v=(v^{s},v^{c},v^{u})\in \mathbb{R}^{s}\oplus\mathbb{R}\oplus\mathbb{R}^{u}=T_{x}M| \| v^{s}+v^{c}\|\leq\alpha\|v^{u}\|\}.$$
  For $\alpha\in(0,1)$, one can check that $\cC_{\alpha}^{s}(x) $ is transverse to $\cC_{\alpha}^{u}(x)$ and $\cC_{\alpha}^{uu}(x)$ is contained in $\cC_{\alpha}^{u}(x)$ for any $x\in \R^n$.

Denote by  $C=[-1,1]^s\times[-1,1]\times[-1,1]^u$.
A \emph{blender horseshoe} $\La$ of $u$-index $i+1$ is a hyperbolic basic set of $u$-index $i+1$ for an embedding map $f: C\mapsto \R^n $  such that:
\begin{itemize}\item[H1)]
the maximal invariant set in $C$ is $\La$ and the dynamics restricted to $C$ is a two-leg horseshoe ( hence it  exhibits  two fixed points $P$ and $Q$), that is, the intersection $f^{-1}(C)\cap C$ consists in two horizontal disjoint sub-cubes $A, B$  and the images $f(A), f(B)$ are two vertical sub-cubes;
     \item[H2)] the set $f(C)\cap [-1,1]^s\times\R\times[-1,1]^u$ consists in two connected components $\cA, \cB$ such that $P\in\cA$ and $Q\in\cB$. Furthermore,
     there exists $\alpha\in(0,1)$ such that  the cone field $\cC_{\alpha}^s$ is strictly $Df^{-1}$   invariant  and the cone fields
     $\cC_{\alpha}^{u}, \cC_{\alpha}^{uu}$ are  strictly $Df$ invariant. Moreover, for  any $x\in   f^{-1}(\cA\cup\cB)$ and   any vector $v\in\cC_{\alpha}^{u}(x)$, $v$ is uniformly expanded by $Df$. Similarly, the vector in $\cC^{s}_{\alpha}$ is uniformly contracted by $Df$.

          \item[H3)] A compact disc $D^{u}$ of dimension $i$ is called \emph{a uu-disc}, if the relative interior of $D^{u}$ is contained in the interior of $C$, $D^u$ is tangent to $\cC_\alpha^{uu}$ and  the boundary of $D^u$ is contained in $[-1,1]^s\times[-1,1]\times\partial{[-1,1]^u}$. Then every uu-disc intersecting $W^s_{loc}(P)$ (resp. $W^s_{loc}(Q)$) is disjoint from $W^s_{loc}(Q)$ (resp. $W^s_{loc}(P)$).
              \item[H4)]   A uu-disc $D^{u}$  is \emph{between $W^s_{loc}(P)$ and $W^s_{loc}(Q)$}, if $D^u$ is homotopic to $W^{uu}_{loc}(P)$ in the set of uu-discs whose homotopy process is disjoint from $W^s_{loc}(Q)$ and $D^u$ is homotopic to $W^{uu}_{loc}(Q)$ in the set of uu-discs whose homotopy process is disjoint from $W^s_{loc}(P)$. Then, for any  uu-disc $D^u$ \emph{ between $W^s_{loc}(P)$ and $W^s_{loc}(Q)$},  at least one of the connected components of $f(D^u)\cap C$
              is a uu-disc between $W^s_{loc}(P)$ and $W^s_{loc}(Q)$
\end{itemize}
The existence of blender horseshoe is a robust property.
The items H1) and H2) imply that there exists a dominated splitting $T_{\La}M=E^s\oplus E^{cu}\oplus E^{uu}$ such that $\dim(E^{uu})=i$ and $\dim(E^{cu})=1$.

  A uu-disc $D^{u}$ is said to be in the \emph{characteristic region}, if $D^{u}$ is between $W^s_{loc}(P)$ and $W^s_{loc}(Q)$. According to item H4) above, for any uu-disc $D^{u}$ in the characteristic region, $f(D^u)$ contains a uu-disc in the characteristic region.

By items H1) and H2), there exists $\epsilon_1>0$ such that
$$\cA\cup\cB\subset (-1+\epsilon_1,1-\epsilon_1)^s\times [-1,1]\times [-1,1]^u$$  and $$f^{-1}(\cA\cup\cB)\subset [-1,1]^s\times [-1,1]\times (-1+\epsilon_1,1-\epsilon_1)^u.$$
  A compact disc  $S$ of dimension $i+1$ is called \emph{a $cu$-strip} if $S$ is tangent to the cone field $\cC_{\alpha}^u$ and is  the image of a $C^1$-embedding map $$\Phi:[-1,1]\times[-1,1]^{u}\mapsto [-1+\epsilon_1,1-\epsilon_1]^s\times [-1,1]\times [-1,1]^u$$
    satisfying that $\Phi(\{t\}\times[-1,1]^{u})$ is a uu-disc, for any $t\in[-1,1]$. The  $cu$-strip $S$ is called in the characteristic region if  for any $t\in[-1,1]$, $\Phi(\{t\}\times[-1,1]^{u})$ is between  $W^s_{loc}(P)$ and $W^s_{loc}(Q)$. The uu-discs  $\Phi(\{-1\}\times[-1,1]^{u})$  and $\Phi(\{1\}\times[-1,1]^{u})$ are called \emph{the vertical boundary components of $S$}.
  For any $cu$-strip $S$, we  define the \emph{central length $\ell^c(S)$ of $S$} as the minimum length of all $C^1$ curves in $S$
   joining the two vertical boundary components of $S$.

 In the rest  of this subsection, we fix $f\in\diff^1(M)$ exhibiting  a blender horseshoe $\La$ corresponding to the cube $C$.
 We fix $\tau>1$ such that for any $x\in C\cap f^{-1}(C)$ and $v\in\cC^u_{\alpha}(x)$, one has
 $\norm{Df(v)}\geq \tau\cdot\norm{v}.$
  \begin{lemma}\label{l.blender horseshoe} The blender horseshoe $\La$ is a dynamically defined blender.
  \end{lemma}
  \proof  
  Let $\mathfrak{D}^{\prime}$ be the set of  uu-discs $D$ satisfying that
   \begin{itemize}
     \item $D$ is in the characteristic  region of the blender;
   \item $D$ is contained in  $[-1+\epsilon_1,1-\epsilon_1]^s\times [-1,1]\times [-1,1]^u$.
 \end{itemize}
 By the item H3) above, there exists $\epsilon_2>0$ such that for any $cu$-strip $S_1$ of central length $2\epsilon_2$ which intersects $W^s_{loc}(P)$ and any $cu$-strip $S_2$ of central length $2\epsilon_2$ which intersects $W^s_{loc}(Q)$, we have that $S_1$ and $S_2$ are disjoint.

Since the existence of blender horseshoe is robust, there exists $\epsilon_3>0$ such that for any diffeomorphism $g$ $\epsilon_3$-close to $f$, the continuation $\La_g$ is a blender horseshoe corresponding to the cube $C$. Let $\epsilon=\min\{\epsilon_2,\epsilon_3\}$.

  For any $\delta>0$ small,   we denote by  $\mathfrak{D}_{\delta}$ the set of uu-discs $D^{u}\in\mathfrak{D^{\prime}}$ such that there exists a
  $cu$-strip $S$ which is  disjoint from  $W^s_{loc}(P)\cup W^s_{loc}(Q)$ and is  defined by
  $$\phi:[-1,1]\times[-1,1]^u\mapsto [-1+\epsilon_1,1-\epsilon_1]^s\times [-1,1]\times [-1,1]^u$$
   satisfying that
   \begin{itemize}
   \item $\phi(\{0\}\times[-1,1]^u)=D^u$,
   \item the central length of the $cu$-strips $\phi([0,1]\times[-1,1]^u)$ and $\phi([-1,0]\times[-1,1]^u)$ are $\delta$;
   \end{itemize}
    One can   check  that $$\emptyset\neq\mathfrak{D}_{\delta}\subsetneq  \mathfrak{D}_{\delta^{\prime}}, \textrm{ for any $\delta^{\prime}<\delta$ small.} $$

  \begin{claim}There exists $\delta_0>0$ small such that $\mathfrak{D}_{\delta_0}$ is a strictly invariant family.
  \end{claim}
  \proof We denote by $c=\sup_{x\in C}\norm{Df(x)}>\tau$. By the uniform expansion  of $Df$ along $\cC^u_{\alpha}$,  there exists $\delta_1<\frac{\epsilon}{2c}$ small enough such that for any $\delta\leq\delta_1$ and
   any uu-disc $D^{u}\in\mathfrak{D}_{\delta}\backslash\mathfrak{D}_{2\delta}$,
we  have that $f(D)$ contains a uu-disc in $\mathfrak{D}_{\tau\cdot\delta}$.

Let $\delta_0=\frac{1}{2} \delta_1$ and $\delta^{\prime}=\tau\delta_0<\epsilon$.
We will prove that for   any $D \in\mathfrak{D}_{2\delta_0}$, $f(D)$ contains a uu-disc in $\mathfrak{D}_{\delta^{\prime}}$.
 Since the blender horseshoe is a horseshoe with two legs,   for any $D\in\mathfrak{D}_{\delta_0}$,
 one has  that $f(D)$ contains two  discs $D_1, D_2$  such that
 \begin{itemize}
 \item $D_1=f(D)\cap\cA$ and $D_2=f(D)\cap \cB$;
 \item $D_1,D_2$ are tangent to the cone field $\cC^{uu}_{\alpha}$.
 \end{itemize}

 By item H4), either $D_1$ or $D_2$ is in the characteristic region of the blender. Without loss of generality, we assume that $D_1$ is in the characteristic region. If $D_1$ is contained in $\mathfrak{D}_{\delta^{\prime}}$, we are done.
  Otherwise,
  one can do a $\delta^{\prime}$-perturbation $\tilde{f}$ of $f$ in $\cA$ without changing $W^s_{loc}(P)$ such that $\tilde{f}(D)\cap \cA$ is not in the superposition region of $\La_{\tilde{f}}$, then one has that $D_2$ is in the superposition of $\La_{\tilde{f}}$ as well as  of $\La$. In this case, one has that $D_2$ must be in $\mathfrak{D}_{\delta^{\prime}}$, otherwise, one can do another  $\delta^{\prime}$-perturbation $f^{\prime}$ of $\tilde{f}$, supported in $\cB$ without changing $W^s_{loc}(Q)$, such that $f^{\prime}(D)\cap\cB$ is not in the superposition region of $\La_{f^{\prime}}$ which contradicts to the item H4) for $f^{\prime}$.

 By  the strictly invariant property of the strong unstable cone field, one has that $\mathfrak{D}_{\delta_0}$ is a strictly invariant family.
  \endproof
Let $\mathfrak{D}$ be the restriction of the family $\mathfrak{D}_{\delta_0}$ to the region $[-1,1]^s\times [-1,1]\times [-1+\epsilon_1,1-\epsilon_1]^u$, this gives a strictly invariant family $\mathfrak{D}$ in the interior of $C$, ending the proof of Lemma~\ref{l.blender horseshoe}.
  \endproof
According to the Lemma \ref{l.blender horseshoe}, 
we  can also denote  a cu-blender horseshoe as  $(\La, C, \cC_{\alpha}^{uu}, \mathfrak{D})$.  $C$ is also called the domain of $\La$.
 Let $\epsilon_0$ be the strength of $\mathfrak{D}$. In the whole paper, for a blender horseshoe, we use the type of strictly invariant family given by Lemma~\ref{l.blender horseshoe}.

Recall that $\tau>1$ is the number such that for any $x\in C\cap f^{-1}(C)$, we have that $$\norm{Df(v)}>\tau\cdot\norm{v}, \textrm{ for any $v\in\cC_{\alpha}^{u}(x)$.}$$
We can prove that  the central length of any cu-strip, ``crossing"  the characteristic  region, is uniformly expanded by the dynamics. To be precise:
 \begin{lemma}\label{l.expanding in the center} For any cu-strip $S$ defined by $\phi:[-1,1]\times[-1,1]^u\mapsto C$ satisfying that $\phi(\{t_0\}\times[-1,1]^u)$ contains an element of $\mathfrak{D}$, for some $t_0\in[-1,1]$,  if the central length of $S$ is less than $\epsilon_0$, one has that $f(S)$ contains a $cu$-strip $S_1$ in $C$
 such that
 \begin{itemize}\item $\ell^c(S_1)>\tau\cdot\ell^c(S)$;
 \item $S_1$ is foliated by discs, one of which contains an element of $\mathfrak{D}$.
 \end{itemize}
\end{lemma}
\proof  By the strictly invariant property of $\mathfrak{D}$,
 $f(\phi(\{t_0\}\times[-1,1]^u))$ contains  a disc $D_1^{u}$ in $\mathfrak{D}$. By the item H2), the connected component of $f(S)\cap C$ containing $D^{u}_1$ is a cu-strip $S_1$ in $C$.  For any $C^1$-curve $\gamma$ in $S_1$ joining the two vertical boundary components of $S_1$, we have that $f^{-1}(\gamma)$ is a $C^1$ curve in $S$ and   joins the two vertical boundary components of $S$. Hence, we have that the length $\ell(f^{-1}(\gamma))$ is no less than $\ell^c(S)$. Since $S$ is tangent to $\cC_{\alpha}^{u}$, by the uniform expansion of $Df$ along the cone field    $\cC_{\alpha}^{u}$, we have that $\ell(\gamma)>\tau\cdot \ell(f^{-1}(\gamma))$, which implies that $\ell^c(S_1)>\tau\cdot\ell^c(S)$.
\endproof
\begin{Remark}  Lemma~\ref{l.expanding in the center} allows us to iterate the cu-strip crossing the superposition region and to gain some expansion in the center direction. This  is also the reason why we use blender horseshoe instead of  the  more general dynamically defined blender.
\end{Remark}

By Lemma \ref{robustness of blender}, we have that for any $D^u\in\mathfrak{D}$, there exists a non-empty intersection of $D^{u}$ and $W^s_{loc}(\La)$.
 Since $\La$ is a hyperbolic basic set, we have that there exists $x\in\La$ such that $D^{u}$ intersects $W^s_{loc}(x)$.  In general, $x$ is not a periodic point.
 The following lemma shows that we can enlarge the disc $D^{u}$ in the center direction, and  after (uniformly) finite many iterations,   the enlarged disc always intersects the local  stable manifolds of periodic orbits in $\La$.
  \begin{lemma}\label{cu strip covers center}
  Let $(\La, C, \cC_{\alpha}^{uu},\mathfrak{D})$ be a cu-blender horseshoe and $\epsilon$ be the strength of $\mathfrak{D}$.
 Then there exists $N\in\N$ such that for any $cu$-strip $D^{cu}$ defined by $\phi:[-1,1]\times[-1,1]^u\mapsto C$ satisfying that
 \begin{itemize}
 \item the central length of $D^{cu}$ is no less than $\epsilon$;
 \item  each $uu$-disc $\phi(\{t\}\times[-1,1]^u)$ contains an element of $\cV_{\epsilon/2}(\mathfrak{D})$;
 \end{itemize}
and for any hyperbolic periodic orbit $\cO_p\subset\La$, one has that
 there exists a point $x\in D^{cu}$ whose forward orbit is in $C$, such that $f^{N}(x)\in W^s_{loc}(\cO_p)$.
  \end{lemma}
  \proof
   By the hyperbolicity of $\La$,   there exist  $\eta>0$ and $\delta>0$ small enough such that for any points $x,y\in\La$ with $\ud(x,y)<\delta$, one has that  $W^s_{\eta}(x)$ is contained in $C$ and intersects $W^u_{\eta}(y)$ transversely in a unique point.

  Since for any $D^u\in\cV_{\epsilon/2}(\mathfrak{D})$, the stable manifold $W^{s}_{loc}(\Lambda)$ intersects  $D^u$ transversely. Hence,  there exists $\delta_0>0$ such that
  for any cu-strip $D^{cu}$ satisfying the hypothesis of Lemma~\ref{cu strip covers center},
   there exists $x\in\La$ such that $W^s_{loc}(x)$ intersects $D^{cu}$ transversely  in a point
    whose distance to the boundary of $D^{cu}$ is no less than $\delta_0$.

    By the uniform continuity of the local stable manifold, there exists $\delta_1<\delta$ such that for any points $x,y\in\La$ with $\ud(x,y)<\delta_1$ and any  disc $D$ with inner  radius  no less than $\delta_0/2$ which is   tangent to the center unstable cone field
     and is centered at a point in $W^{s}_{loc}(x)$,
     one has  that $W^s_{loc}(y)$ intersects $D^{cu}$ transversely in the interior.

    Since $\La$ is a hyperbolic basic set, there exists a periodic orbit $\cO_{p_0}$ and a positive integer $N$ such that
    \begin{itemize}
    \item the orbit $\cO_{p_0}$  is $\delta_1/2$ dense in $\Lambda$.
      \item  for any two points $p_1,p_2\in\cO_{p_0}$, there exists an integer $n\in[0,N]$
    such that $f^{-n}(W_{\eta}^u(p_1))\subset W^u_{\delta_1/2}(p_2)$.
    \end{itemize}

    By the choice of $N$,  one has that $N$  depends only on the number $\delta_0$ and on the set $\La$ and one can check that
    $N$ is the integer that we need,  ending the proof of Lemma \ref{cu strip covers center}.
\endproof

\subsection{Flip-flop configuration}

    In this section, we recall the definition and the properties of  \emph{flip-flop configuration}. Roughly speaking,  a flip flop configuration   is  a robust cycle formed by a cu-blender and a hyperbolic periodic orbit of different indices in the way:  the unstable manifold of the periodic orbit ``crosses" the superposition region of cu-blender, and every disc in the strictly invariant family intersects  the stable manifold of the periodic orbit. To be specific:
    \begin{definition}\label{d.f} Consider a dynamically defined blender $(\Lambda,U,\mathcal{C}^{uu},\mathfrak{D})$  of uu-index $i$ and  a hyperbolic periodic point  $q$  of u-index $i$. We say that $\Lambda$ and $q$ form a \emph{flip-flop configuration}, if there exist a disc $\Delta^{u}\subset W^{u}(q)$ and a compact  submanifold with boundary $\Delta^{s}\subset W^{s}(q)\cap U$ such that:
    \begin{enumerate}\item $\Delta^{u}\in\mathfrak{D}$ and $f^{-n}(\Delta^{u})\cap \overline{U}=\emptyset$, for any $n\in\N^+$;
    \item there exists $N\in\N$ such that for any integer $n>N$, $f^{n}(\Delta^{s})\cap\overline{U}=\emptyset$; Moreover, for any $x\in\Delta^s$, if $f^{j}(x)\notin\overline{U}$ for some $j>0$,  the forward orbit of $f^{j}(x)$ is in the complement  of $\overline{U}$;
        \item for any $y\in\Delta^{s}$, $T_yW^s(q)\cap\mathcal{C}^{uu}=\{0\}$;
        \item there exist a compact set $K\subset\Delta^{s}$ and a  number $\eta>0$ such that for any disc $D\in\mathfrak{D}$, the disc $D$ intersects $K$ in  a point whose distance to $\partial{D}$ is no less than $\eta$.
    \end{enumerate}
We denote the flip flop configuration as $(\La, U, \mathcal{C}^{uu},\mathfrak{D},\cO_q, \Delta^s,\Delta^u)$.
    \end{definition}
    \begin{Remark}
     \begin{enumerate}
     \item By the  item (4) in Definition~\ref{d.f}, for the discs in $\mathfrak{D}$,  their diameters are uniformly  bounded away from zero;
    \item  It's shown in \cite[Proposition 4.2]{BBD2} that the existence of flip-flop configuration is a robust property.
    \end{enumerate}
    \end{Remark}

    One says that a set $V$ is a \emph{neighborhood of the  flip-flop configuration \\ $(\La, U, \mathcal{C}^{uu},\mathfrak{D},\cO_q, \Delta^s,\Delta^u)$},  if its interior contains the set  $$\mathcal{O}_q\cup\overline{U}\cup\,\bigcup_{j\geq 0}f^{j}(\Delta^s)\cup\,\bigcup_{j\geq 0}f^{-j}(\Delta^u).$$

    \begin{lemma}\label{l.strong unstable cone field}
    \cite[Lemma 4.6] {BBD2}Let $f\in\diff^1(M)$. Assume that there exists a  flip-flop configuration
    $(\La, U,\mathcal{C}^{uu},\mathfrak{D},\cO_q, \Delta^s,\Delta^u)$. For any small enough compact neighborhood $V$ of  the flip-flop configuration, one has that  the maximal invariant set of $V$ has a partially hyperbolic splitting of the form $E^{cs}\oplus E^{uu}$, where $\dim(E^{uu})$ equals the u-index of $q$. Moreover, there exists a strictly $Df$-invariant cone field $\mathcal{C}_{V}^{uu}$ over $V$ which continuously extends the cone field $\mathcal{C}^{uu}$,   such that any vector in $\mathcal{C}_{V}^{uu}$ is uniformly expanded by $Df$.
    \end{lemma}
    \begin{definition}\label{d.split flip flop}
    Given  a flip-flop configuration  $(\La, U,\mathcal{C}^{uu},\mathfrak{D},\cO_q, \Delta^s,\Delta^u)$. Let $i$ be the u-index of the periodic point $q$.
     We say that this configuration is \emph{split} if there exists a small compact neighborhood $V$ of this configuration such that the maximal invariant set of $V$ admits a partially hyperbolic splitting of the form $E^{ss}\oplus E^{c}\oplus E^{uu}$, where $\dim(E^{ss})=ind(\Lambda)$ and $\dim(E^{uu})=i$.
    \end{definition}
      The following proposition gives the existence of split flip-flop configuration, whose proof can be found  in \cite[Section 5.2]{BBD2}.
    \begin{proposition}\label{p.existence of flip flop} \cite{BBD2}
    Let $\cU$ be an open set of diffeomorphisms such that for any $f\in\cU$, there exist two hyperbolic periodic points $p_f, q_f$ of  u-index $i_p>i_q$ respectively which continuously depend on $f$ and are in the same chain class $C(p_f, f)$.

    Then there exists an open and dense subset $\tilde{\cU}$ of $\cU$ such that for any $f\in\tilde{\cU}$ and any $i\in(i_q,i_p]$,  there exists a split flop-flop configuration formed by a dynamically defined cu-blender of uu-index $i-1$ and a hyperbolic periodic orbit of u-index $i-1$.
    \end{proposition} 
\section{Approximation of non hyperbolic ergodic measure  by periodic orbits: Proof of Theorem~\ref{thmF}.}
    In this section, we prove that  all the non-hyperbolic ergodic measures supported  in a small enough neighborhood of the split flip-flop configuration are approximated by hyperbolic periodic measures.

Consider $f\in\diff^1(M)$.  Assume that there exists a split flip flop configuration formed by  a cu-blender horseshoe $(\La,C,\cC^{uu},\mathfrak{D})$ of u-index $i+1$ and a hyperbolic periodic point $q$ of u-index $i$.
 Let $\epsilon_0$ be the strength of the strictly invariant family $\mathfrak{D}$.
We denote by   $\cC^{u}$  the center  unstable invariant cone field   defined in $C$.
By  the definition of blender horseshoe,   there exists $\tau_1>1$ such that
$$\norm{Df(v)}>\tau_1\cdot \norm{v}, \textrm{ for any $x\in C\cap f^{-1}(C)$ and $v\in\cC^{u}(x)$.}$$
 By Lemma~\ref{l.strong unstable cone field} and the definition of split flip flop configuration, there exist a small neighborhood  $V$ of the configuration such that
 \begin{itemize}
 \item the maximal invariant $\tilde\La$ in $\overline{V}$ admits a partially hyperbolic splitting $T_{\tilde\La}M=E^{ss}\oplus E^c\oplus E^{uu}$, where
  $\dim (E^c)=1$ and $\dim (E^{uu})=i$.
 \item there exist a continuous extension $\cC^{uu}_V$ of $\cC^{uu}$ in $V$ and a number $\tau_2>1$ such that
$$\norm{Df(v)}\geq \tau_2\cdot\norm{v}, \textrm{ for any $x\in V\cap f^{-1}(V)$ and any  $v\in \cC^{uu}(x)$}.$$
\end{itemize}
We assume, in addition, that there exists a $Df$-strictly invariant cone field $\cC^{u}_V$ in $V$ which is a continuous extension of the center unstable cone field $\cC^u$ in $C$.

We denote by $b=\max\{\sup_{x\in M}\norm{Df_x}, \sup_{x\in M}\norm{Df^{-1}_x}\}.$
  \begin{proposition}\label{p.main} With the notation above.
  Let $\nu$ be a non-hyperbolic ergodic measure supported on $\tilde{\La}$. Assume that there exists a periodic orbit $\mathcal{O}_p\subset \La$ and a point $y\in\tilde{\La}\cap W^u_{loc}(\cO_p)$ in the basin of $\nu$ such that $\overline{Orb(y,f)}$ is far away from the boundary of $V$.

   Then $\nu$ is approximated by hyperbolic periodic orbits which are homoclinically related to $\La.$
  \end{proposition}

  \proof  

By plaque family theorem, we fix plaque families $W^{cu}$ and $W^{cs}$  for $\tilde\La$
such that $W^{cu}(y)$ is foliated by  discs tangent to cone field $\mathcal{C}_V^{uu}$.
 Since $y\in W^{u}_{loc}(\cO_p)$ and $\Lambda$ is uniformly hyperbolic,   there exists an integer $N_0\in\N$, which is only depends on the  size
 of $W^u_{loc}(\cO_p)$ and $\Lambda$,
 such that $f^{-N_0}(y)$ is a $(\tau_1, E^c\oplus E^u)$ expanding point. By Lemma \ref{lp}, there exists $\delta_0$ independent of $y$ such that
 $W^{cu}_{\delta_0}(y)$ is contained in the unstable manifold of $y$, which implies that $W^{cu}_{\delta_0}(y)\subset W^{u}_{loc}(\cO_p)$.

 For any point $z\in\tilde\La$, we denote by
  $$W^{c}(z)=W^{cu}(z)\cap W^{cs}(z),$$
 then this  gives a plaque family for $E^c$. For any $z\in\tilde{\La}$ and any point $w\in W^i(z)$, we denote by $\tilde{E}^i(w,z)=T_w W^i(z)$ for $i=c,cu$.
 When there is no ambiguity,  for $i=c,cu$, we denote $\tilde{E}^i(w,z)$ as $\tilde{E}^i(w)$ for simplicity.

To fulfill the proof, we need the following lemma:
\begin{lemma}\label{l.orbit follow the measure} For any $\varepsilon>0$,  there exist an integer $N>0$ and  a sequence of points $\{z_k\}_{k\in \N}\subset W^u_{loc}(\cO_p)$ together with a sequence of  positive integers $\{t_k\}_{k\in\N}$   such that
 \begin{itemize}
 \item
 $$f^{t_k }(z_k)\in W^s_{loc}(\cO_p)\textrm{ and } \ud\big(\frac{1}{t_k}\sum_{j=0}^{t_k-1}\delta_{f^j(z_k)},\nu\big)<\varepsilon;$$
 \item For any $j\in[N, t_k]$, we have that $$-\varepsilon<\frac{1}{j}\sum_{i=0}^{j-1}\log\norm{Df|_{_{E^c(f^j(z_k))}}}<\varepsilon;$$
 \item the sequence $\{t_k\}_{k\in\N}$ tends to infinity and the orbit segments $\big\{\{z_k, t_k\}\big\}_{k\in\N}$ are contained in $V$.
 \end{itemize}
\end{lemma}
\proof
For any  $\epsilon>0$ small,   by the uniform continuity of $Df$ on the unit tangent bundle of $TM$, and the compactness of $M$ and $\tilde{\La}$,
 there exists $\delta\in(0,\delta_0)$ such that:
\begin{itemize}
\item[--] Given $x\in\tilde\La$. For any  two points $x_1,x_2\in W^c_{\delta}(x)$,   one has that
              $$ -\epsilon \leqslant \log{\norm{Df|_{\tilde{E}^{c}(x_1)}}}-\log{\norm{Df|_{\tilde{E}^{c}(x_2)}}}\leqslant  \epsilon;$$
\item[--] for any point $z\in\tilde\La$ and any point $w\in W^{cu}_{\delta}(z)$, one has that
 $$ -\epsilon\leq\log{\um(Df|_{\tilde{E}^{cu}(w)})}-\log{\norm{Df|_{E^c(z)}}}\leq\epsilon;$$
\item[--] for any $w_1,w_2\in\tilde{\La}$ satisfying that $\ud(w_1,w_2)<\delta$, we have that $$|\log\norm{Df|_{_{E^c(w_1)}}}-\log\norm{Df|_{_{E^c(w_2)}}}|<\epsilon/2;$$
\item[--] for any two points $z_1,z_2\in M$ satisfying that $\ud(z_1,z_2)<\delta$, we have that $\ud(\delta_{z_1},\delta_{z_2})<\epsilon/2,$
                 where $\delta_{z_i}$ denotes the Dirac measure supported on the point $z_i$;
   \item[--] for any point $x\in\tilde\La$, one has that $f(W^{cu}_{\delta}(x))\subset W^{cu}(f(x))$.
\end{itemize}
   By the choice of $y$,  there exists  an integer $N$ such that for any $n\geqslant N$, we have that
    \begin{equation}\label{equ:estimate on y}-\epsilon<\frac{1}{n}\sum_{i=0}^{n-1}\log\norm{Df|_{E^{c}(f^{i}(y))}}<\epsilon \textrm{ and } \ud \big(\,\frac{1}{n}\sum_{i=0}^{n-1}\delta_{f^{i}(y)},\nu\,\big)<\epsilon/2.
    \end{equation}



     For any $C^1$ curve $\gamma\subset M$, we denote by $\ell(\gamma)$ the length of $\gamma$.
For any $n\in\N$, we take the  $C^1$ curve $\gamma_{n}\subset W^{c}(y)$ such that
     \begin{itemize}
     \item $\ell(\gamma_{n})=\delta\cdot e^{-4n\cdot\epsilon}$;
     \item  the curve  $\gamma_n$ is centered at $y$.
     \end{itemize}

     \begin{claim}\label{c.length of central curve}  There exists an integer $N_1$ such that for any $n>N_1$, we have that
     $$\ell(f^{i}(\gamma_{n}))<\delta,\textrm{ for any integer $i\in[0,n]$. }$$
    \end{claim}
    \proof   Recall that $b\geq\max_{x\in M}\norm{Df_x}$, then
    there exists an integer $N_1$ satisfying  that
    $$e^{-n\epsilon}\cdot b^{N}<1, \textrm{ for any $n>N_1$,}$$
Hence, for any  $n>N_1$ and any integer $i\in[0,N]$, we have the estimate:
    $$\ell(f^{i}(\gamma_{n}))=\int_{0}^{1}\|\frac{\ud}{\ud t}f^{i}(\gamma_{n}(t))\|\ud t  \leq  b^{i}\cdot \ell(\gamma_{n})<\delta;$$

    We will   prove this claim  for $i\in(N, n]$ inductively.
      Assume that for any integer $j\leq i\in(N,n)$, we have that $\ell(f^{j}(\gamma_{n}))<\delta$,
      then by the choice of $\delta$ and $N$, we have that
      \begin{align*}
       \ell(f^{i+1}(\gamma_{n}))=\int_{0}^{1}\|\frac{\ud}{\ud t}f^{i+1}(\gamma_{n}(t))\|\ud t
       &\leq\int_{0}^{1}\prod_{j=0}^{i}\norm{ Df|_{\tilde{E}^{c}(f^{j}(\gamma_n(t)))} } \|\gamma_{n}^{\prime}(t)\|\ud t
       \\
         &\leq\int_{0}^{1}e^{(i+1)\epsilon}\cdot\prod_{j=0}^{i}\norm{ Df|_{E^{c}(f^{j}(y))} } \|\gamma_{n}^{\prime}(t)\|\ud t
         \\
         &<\delta.
         \end{align*}
       \endproof
    By  Claim \ref{c.length of central curve} and the choice of $\delta$, for any $n>N_1$,  we have the estimate:
\begin{align*}\ell(f^{n}(\gamma_{n}))&=\int_{0}^{1}\|\frac{\ud}{\ud t}f^{n}(\gamma_{n}(t))\|dt
\\
 &\geq\int_{0}^{1}e^{-n \epsilon}\prod_{j=0}^{n-1}\norm{ Df|_{E^{c}(f^{j}(y))} } \|\gamma_{n}^{\prime}(t)\|\ud t
 \\
 &
 \geq \delta \cdot e^{-6n\epsilon}.
 \end{align*}

A $i+1$-dimensional disc $D^{cu}\subset V$  is called a \emph{ uu-foliated cu-disc}, if  one has that
\begin{itemize}
\item the disc $D^{cu}$ is tangent to the center unstable cone field $\cC_V^{u}$
\item there exists a $C^1$
  embedding $\phi:[0,1]\times[0,1]^i\mapsto M$ such that $\phi([0,1]\times[0,1]^i)=D^{cu}$ and
  for any $t\in[0,1]$,  $\phi(\{t\}\times [0,1]^{i})$ is a disc tangent to  $\cC_V^{uu}$.
  \end{itemize}
  The \emph{ central length $\ell^c(D^{cu})$ of   $D^{cu}$} is defined as the infimum of the  length of the $C^1$ curves contained in $D^{cu} $ joining the discs
   $\phi(\{0\}\times [0,1]^{i})$ and $\phi(\{1\}\times [0,1]^{i})$. Two discs $\phi(\{0\}\times [0,1]^{i})$ and $\phi(\{1\}\times [0,1]^{i})$ are called \emph{the vertical boundary components of $D^{cu}$}.

     Consider a submanifold
    $S_{n}$ which is the   $\ell(\gamma_n)$ tubular neighborhood of $W^{uu}_{\delta}(y)$  in $W^{cu}(y)$, then  $S_n$ is a uu-foliated cu-disc for $n$ large and is contained in $W^{u}_{loc}(\cO_p)$, moreover, the central length of $S_n$ is   $\delta \cdot e^{-4n\epsilon}.$ We denote by $S_n(i)$   the connected component of $f^{i}(S_n)\cap B_{\delta}(f^{i}(y))$ which contains $f^{i}(y)$ (see Figure 1).

      \begin{figure}[h]
\begin{center}
\def\svgwidth{0.8\columnwidth}
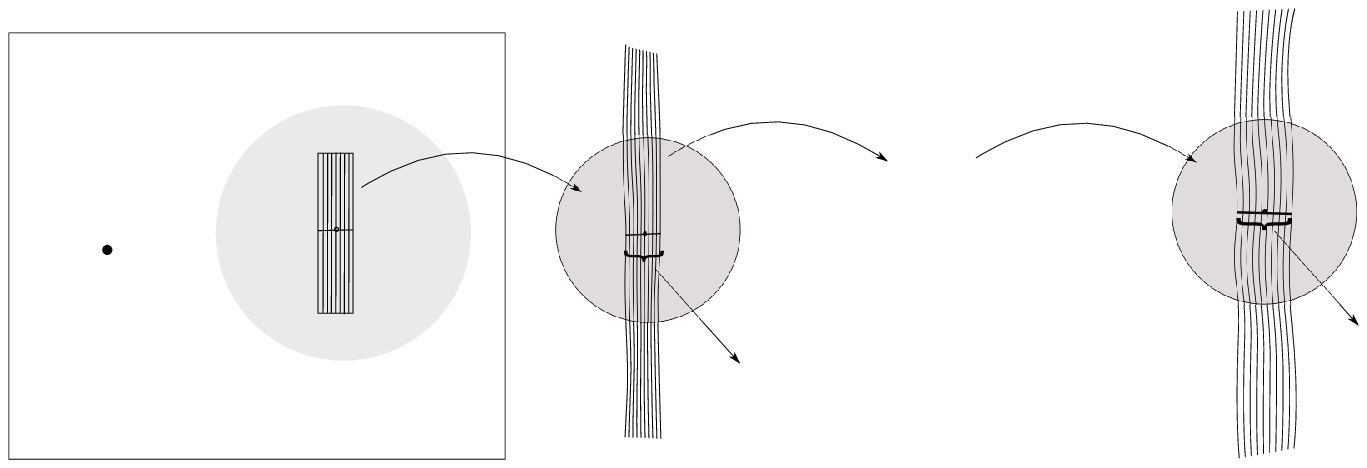
  \caption{}
\end{center}
\end{figure}

    By Claim \ref{c.length of central curve} and uniform expansion of $Df$ along the cone field $\cC_V^{uu}$, one has that
    \begin{itemize}
    \item $S_n(i)$ is a uu-foliated  cu-disc whose vertical boundary components  are contained in the $f^i$-image of the vertical boundary components of $S_n$;
    \item $S_n(i)$ is saturated by the discs tangent to $\cC^{uu}_V$ with diameter of size  $\delta$, for $i=1,\cdots,n$.
      \end{itemize}
    \begin{claim}\label{c.lower bound of central length} The central length of $S_n(i)$ is no less than $\delta\cdot e^{-6n\epsilon}$, for each $i\in[1, n]$.
    \end{claim}
    \proof Given $i\in[1,n]$, for any $C^1$ curve $\xi(t)_{t\in[0,1]}$ contained in $S_n(i)$ which joins the two vertical boundary components  of $S_n(i)$,  one has  that $f^{-i}(\xi(t))$ is a $C^1$ curve joining the vertical  boundary components  of $S_n$.
    Moreover, by the definition of $S_n(i)$, one has that for any $j\in[1,i]$, one has that $f^{-i+j}(\xi(t))\subset W^{cu}_{\delta}(f^j(y))$.

    Since $\ell(f^{-i}(\xi(t)))\geq \delta\cdot e^{-4n\epsilon}$, by the choice of $\delta$, one has  that
    \begin{align*}
    \ell(\xi(t))&=\ell(  f^{i}\circ f^{-i}(\xi(t)))
    \\
    &=\int_0^{1} \norm{  Df^{i}\frac{\ud}{\ud t}f^{-i}(\xi(t))}\ud t
    \\
    &\geq \int_0^1e^{-n\epsilon}\cdot \prod_{j=0}^{i-1}\norm{Df|_{E^c(f^j(y))}}\norm{\frac{\ud}{\ud t}f^{-i}(\xi(t))}\ud t
    \\
    &\geq e^{-2n\epsilon}\cdot\ell(f^{-i}(\xi(t)))
    \\
    &\geq \delta\cdot e^{-6n\epsilon}
 \end{align*}
 \endproof

  By the definition  of   flip-flop configuration and the choice of $y$,
there exists a  sequence of positive integers $\{n_{k}\}_{k\in\N}$ tending to infinity,
such that $f^{n_{k}}(y)$ is in a small neighborhood of $\Delta^{u}$, hence $f^{n_k}(y)$ is in $ C$.
Consider the $cu$-disc  $S_{n_{k}}$, by  Claim \ref{c.lower bound of central length},
    we have that $S_{n_k}(n_k)$ is a cu-disc of central length at least $\delta\cdot e^{-6\cdot n_k\cdot \epsilon}$.

Recall that $\tau_2>1$ is a number such that for any $x\in V\cap f^{-1}(V)$ and $v\in\cC^{uu}_V(x)$, one has that $\norm{Df(v)}\geq\tau_2\norm{v}$. We  denote by $N(\delta)$ the smallest integer satisfying that $\tau_2^{N(\delta)}\cdot \delta\geq b_0$, where $b_0$ is the upper bound for the diameters of the uu-discs in
the family $\mathfrak{D}$,  then $N(\delta)\leq \frac{\log{b_0}-\log\delta}{\log\tau_2}+1$.
 By the $Df$-strictly invariant property  of the cone fields $\cC^u_V,\cC^{uu}_V$ and the fact that $\overline{Orb(y,f)}$ is far away from the boundary of $V$,
there exists an integer $\tilde{n}_{k}$ such that
    \begin{itemize}\item $n_k-\tilde{n}_{k}\in(0, N(\delta)]$;
    \item Denote by  $\tilde{S}_{n_k}(n_k)$  the connected component of $f^{n_k-\tilde{n}_{k}}(S_{n_k}(\tilde{n}_{k}))\cap V$ containing
    $f^{n_k}(y)$.
    Then $\tilde{S}_{n_k}(n_{k})$ contains  a cu-strip $\tilde{S}_k$ in the characteristic region of blender horseshoe (see Figure 2) such that for the central length $\ell^c(\tilde{S}_k)$ of $\tilde{S}_k$ , one has
             $$\ell^c(\tilde{S}_k)\geq b^{-N(\delta)}\cdot \delta\cdot e^{-6\cdot n_k\cdot\epsilon}.$$
    \end{itemize}
  \begin{figure}[h]
\begin{center}
\def\svgwidth{0.9\columnwidth}
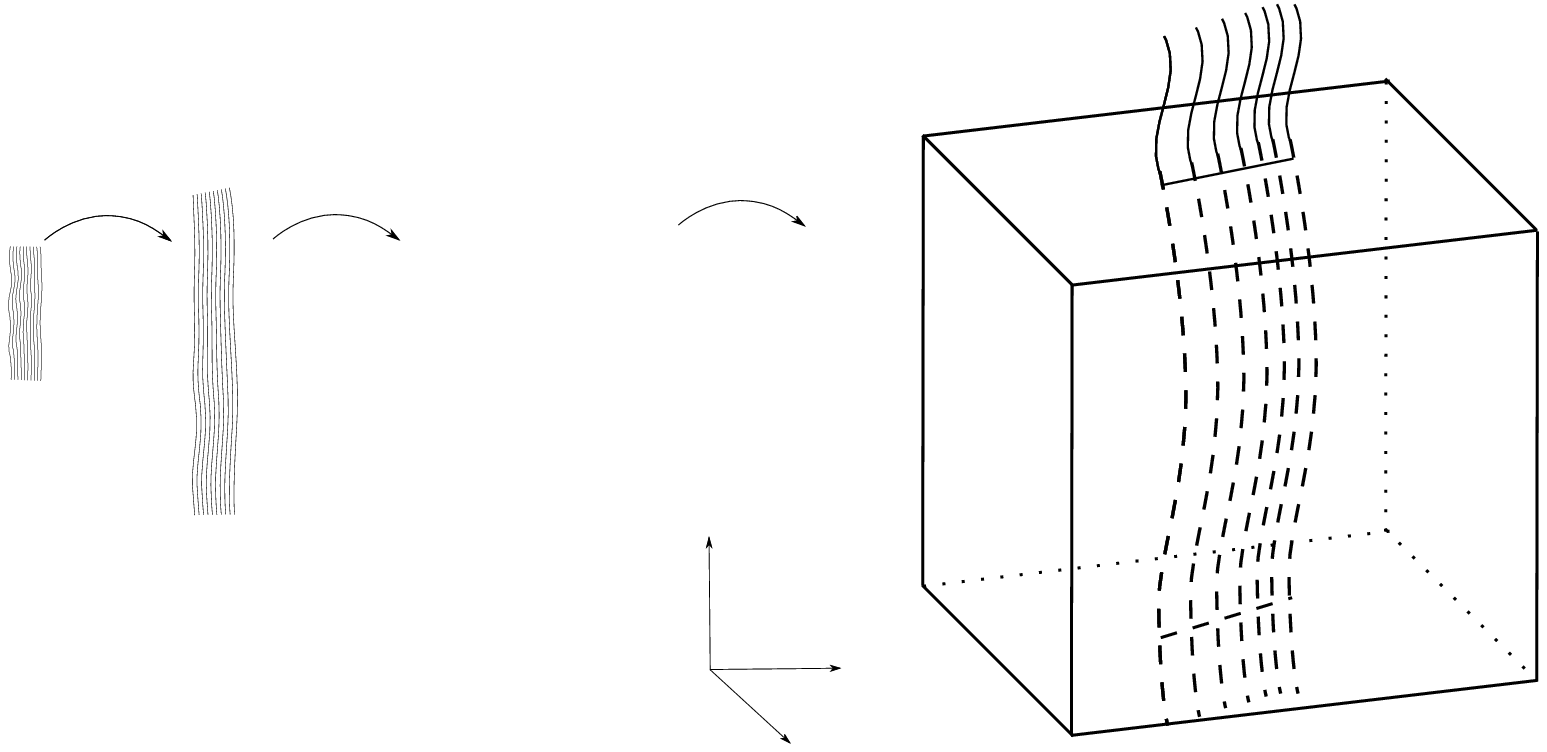
  \caption{}
\end{center}
\end{figure}

   Due to the uniform expansion of $Df$ along the cone field $\cC^{uu}_V$, each point in $f^{-n_k}(\tilde{S}_{k})$
     would stay close to the orbits segment $\{y,\ldots,f^{n_{k}}(y)\}$ for a large proportion
      of time in  $[\, 0,n_{k} \,]\cap\N$; moreover, the proportion  would tend to 1 when $n_{k}$ tends to infinity.

\vspace{2mm}
     We will iterate $\tilde{S}_k$ to make it cut the local stable manifold of $\mathcal{O}_{p}$ transversely.
    Let $\epsilon_0$ be the strength of the strictly invariant family $\mathfrak{D}$,  and we denote by
    $$[r]=\sup\{n\in\Z: n\leq r\}\textrm{ and } N_{k}=\Big[\frac{6\cdot n_{k}\cdot\epsilon}{\log\tau_1}+\frac{\log\epsilon_0-\log\delta+\log{b}\cdot N(\delta)}{\log\tau_1}\Big]+1.$$

   Since the integer $N$ in Lemma \ref{cu strip covers center} is a   constant, for simplicity, we  take its value as $0$.
     By Lemma \ref{l.expanding in the center},  Lemma \ref{cu strip covers  center}
      and the choice of $N_{k}$, there exists an integer $\tilde{N}_{k}\leq N_k$ such that $f^{\tilde{N}_{k}}(\tilde{S}_{k})\cap C$ contains a connected component which intersects the  local stable manifold of $\mathcal{O}_{p}$ in  a point $z$.  Denote by
    $$t_{k}=\tilde{N}_{k}+n_{k} \textrm{ and } z_{k}=f^{-t_{k}}(z),$$
  then one has that  $z_k\in \tilde{\La}$ and satisfies the third item of Lemma \ref{l.orbit follow the measure}.

      Since the choice of $\delta$ is independent of $n_k$, by the fact that $|\tilde{n}_k-n_k|\leq N(\delta)$, we can take $n_{k}$  large enough such that
       $$\frac{N(\delta)+N_{k}}{\tilde{n}_{k}}<\frac{7\epsilon}{\log\tau_1}.$$

By the choice of $\delta$ and $z_k$,  for any $j\in[N, \tilde{n}_k]$, we have that
   \begin{itemize}\item $$\ud(\frac{1}{j}\sum_{i=0}^{j-1}\delta_{f^{i}(z_k)}, \frac{1}{j}\sum_{i=0}^{j-1}\delta_{f^{i}(y)})<\frac{\epsilon}{2};$$
   \item $$\big|\frac{1}{j}\sum_{i=0}^{j-1}\log\norm{Df|_{_{E^c(f^{i}(z_k))}}}- \frac{1}{j}\sum_{i=0}^{j-1}\log\norm{Df|_{_{E^c(f^{i}(y))}}}\big|<\frac{\epsilon}{2}.$$
\end{itemize}

On the other hand, when $n_k$ is large enough, we have that
\begin{align*}
&\ud\big(\frac{1}{t_k}\sum_{i=0}^{t_k-1}\delta_{f^{i}(z_k)}, \frac{1}{\tilde{n}_k}\sum_{i=0}^{\tilde{n}_k-1}\delta_{f^{i}(y)}\big)
\\
&\leq \ud\big(\frac{1}{t_k}\sum_{i=0}^{t_k-1}\delta_{f^{i}(z_k)}, \frac{1}{\tilde{n}_k}\sum_{i=0}^{\tilde{n}_k-1}\delta_{f^{i}(z_k)}\big)
+\ud\big(\frac{1}{\tilde{n}_k}\sum_{i=0}^{\tilde{n}_k-1}\delta_{f^{i}(z_k)}, \frac{1}{\tilde{n}_k}\sum_{i=0}^{\tilde{n}_k-1}\delta_{f^{i}(y)}\big)
\\
&\leq \frac{t_k-\tilde{n}_k}{t_k}\cdot (1+\frac{1}{\tilde{n}_k})+\frac{\epsilon}{2}
\\
&<\frac{7}{\log\tau_1}\epsilon+\frac{\epsilon}{2}
\end{align*}
and for any $j\in(\tilde{n}_k, t_k]$
\begin{align*}&\Big|\frac{1}{j}\sum_{i=0}^{j-1}\log\norm{Df|_{_{E^c(f^{i}(z_k))}}}- \frac{1}{\tilde{n}_k}\sum_{i=0}^{\tilde{n}_k-1}\log\norm{Df|_{_{E^c(f^{i}(y))}}}\Big|
\\
&\leq \Big|\frac{1}{j}\sum_{i=0}^{\tilde{n}_k-1}\log\norm{Df|_{_{E^c(f^{i}(z_k))}}}- \frac{1}{\tilde{n}_k}\sum_{i=0}^{\tilde{n}_k-1}\log\norm{Df|_{_{E^c(f^{i}(y))}}}\Big|
+\Big|\frac{1}{j}\sum_{i=\tilde{n}_k}^{j-1}\log\norm{Df|_{_{E^c(f^{i}(z_k))}}}\Big|
\\
&<\frac{\epsilon}{2}+\frac{7}{\log\tau_1}\cdot \epsilon\cdot \log{b}.
  \end{align*}
As a consequence, we have that
$$\ud(\frac{1}{t_k}\sum_{i=0}^{t_k-1}\delta_{f^{i}(z_k)}, \nu)<\epsilon+\frac{7}{\log\tau_1}\epsilon $$
and
$$\big|\frac{1}{j}\sum_{i=0}^{j-1}\log\norm{Df|_{_{E^c(f^{i}(z_k))}}}\big|<\epsilon+\frac{7\log{b}}{\log\tau_1}\cdot \epsilon, \textrm{for any $j=N,\cdots,t_k.$}$$
Let $c=\max\{1+\frac{7\log{b}}{\log\tau_1},1+\frac{7}{\log\tau_1}\}$, then we only need to   take $\epsilon$ small such that $c\cdot\epsilon<\varepsilon$, ending the proof of Lemma \ref{l.orbit follow the measure}.
\endproof

 Consider the convex sum $\{\alpha\delta_{\mathcal{O}_{p}}+(1-\alpha)\nu\}_{\alpha\in[0,1]}$.
  We fix  $\alpha\in(0,1]$, then  the mean center Lyapunov exponent of $\alpha\delta_{\mathcal{O}_{p}}+(1-\alpha)\nu$ is:
 $$\lambda^{c}(\alpha\delta_{\mathcal{O}_{p}}+(1-\alpha)\nu)=\alpha\lambda^c(\delta_{\cO_p})>0.$$
  Denote by$$\lambda=\exp(-\lambda^{c}(\alpha\delta_{\mathcal{O}_{p}}+(1-\alpha)\nu))\in(0,1).$$

  We take $\varepsilon<\frac{-\log\lambda}{16}$ small. 
   By the uniform continuity of $\log\norm{Df|_{E^c}}$ over $\tilde{\La}$, there exists $\delta_1>0$ such that for any $w_1,w_2\in\tilde{\La}$ satisfying that $\ud(w_1,w_2)<\delta_1$, we have that
  $$\big|\log\norm{Df|_{E^c(w_1)}}-\log\norm{Df|_{E^c(w_2)}}\big|< \frac{-\log\lambda}{16}.$$

  By Lemma \ref{shadow}, we get two numbers $L>0$ and $d_{0}>0$ such that for any $d\in(0,d_0)$, one has that every $\sqrt{\lambda}$-quasi hyperbolic periodic $d$-pseudo orbit   corresponding  to the splitting
  $T_{\tilde\La}M=E^{ss}\oplus (E^c\oplus E^{uu})$ is $L\cdot d$ shadowed by a periodic orbit.

  For any  $d\in(0, \min\{d_{0},\delta_1\})$ small enough   whose  precise value  would be fixed later,   there exists an integer $N_{d}$ such that
 $$f^{-N_{d}}(W^{u}_{loc}(\mathcal{O}_{p}))\subset W^u_{d/2}(\mathcal{O}_{p}) \textrm{  and }
 f^{N_{d}}(W^{s}_{loc}(\mathcal{O}_{p}))\subset W^s_{d/2}(\mathcal{O}_{p}).$$

 By Lemma \ref{l.orbit follow the measure}, there exist an integer $N$, a sequence of points  $\{z_k\}_{k\in\N}$ in $\tilde\La$ and a sequence of integers $\{t_k\}_{k\in\N}$ tending to infinity such that
\begin{itemize}\item $$\ud\big(\frac{1}{t_k}\sum_{j=0}^{t_k-1}\delta_{f^j(z_k)},\nu\big)<  \varepsilon;$$
\item  $$z_k\in W^u_{loc}(\cO_p)\textrm{ and }f^{t_k}(z_k)\in W^s_{loc}(\cO_p);$$

 \item For any $j\in[N, t_k]$, we have that
 \begin{equation}\label{equ:average on center}
 -\varepsilon<\frac{1}{j}\sum_{i=0}^{j-1}\log\norm{Df|_{_{E^c(f^j(z_k))}}}<\varepsilon.
 \end{equation}
 \end{itemize}
     \begin{claim}\label{c.choice of  t_k} there exist integers $m_{k}$ and $t_k$ arbitrarily large such that
      \begin{itemize}
      \item $$|\frac{m_{k}}{t_{k}}-\frac{\alpha}{1-\alpha}|+\frac{2N_d}{t_k}<\varepsilon;$$
      \item The orbit segment  $\{f^{-N_d}(z_{k}),\ldots,f^{N_{d}+t_k+m_{k}}(z_k)\}$ is a $\sqrt{\lambda}$-quasi hyperbolic string corresponding to the splitting
                  $E^{ss}\oplus (E^c\oplus E^{uu})$.
      \end{itemize}
     \end{claim}
  \begin{figure}[h]
\begin{center}
\def\svgwidth{0.8\columnwidth}
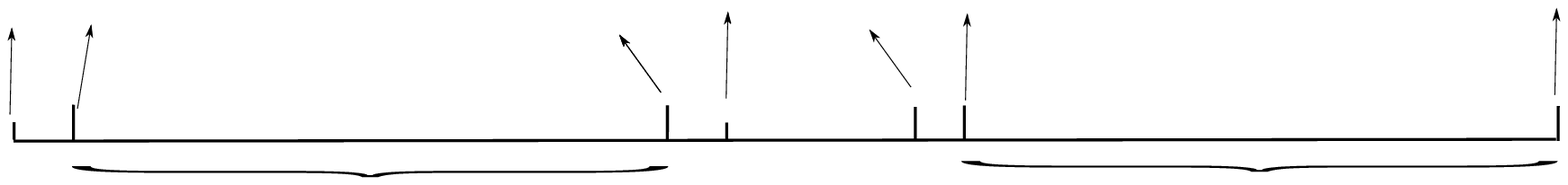
  \caption{}
\end{center}
\end{figure}

\proof

 Since the sequence $\{t_k\}$ tends to infinity, one can take a sequence of positive integers $\{m_k\}_{k\in\N}$ tending to infinity  such that
 $$\lim_{k\rightarrow\infty}|\frac{m_{k}}{t_{k}}-\frac{\alpha}{1-\alpha}|+\frac{2N_d}{t_k}=0,$$
 hence, for $k$ large enough, the first item is satisfied.

 We denote by $\pi_{k,d}=2N_d+t_k+m_k.$
By the choice of $N_d$, we have that $$\ud(f^{-N_d}(z_k),f^{N_{d}+t_k+m_{k}}(z_k))<d.$$
Since $d$ is less than $\delta_1$ and the integers $t_k,m_k$ can be chosen arbitrarily large, by the choice of $\varepsilon$ and Equation (\ref{equ:average on center}),  one has the following estimate:
\begin{align*}
\frac{1}{\pi_{k,d}}\sum_{j=0}^{\pi_{k,d}-1}\log\norm{Df|_{E^c(f^{j-N_d}(z_k))}}
&=\frac{1}{\pi_{k,d}}\sum_{j=0}^{t_k+N_d-1}\log\norm{Df|_{E^c(f^{j-N_d}(z_k))}}
\\
&\hspace{5mm}+\frac{1}{\pi_{k,d}}\sum_{j=t_k+N_d}^{\pi_{k,d}-1}\log\norm{Df|_{E^c(f^{j-N_d}(z_k))}}
\\
&\geq \frac{t_k+N_d}{\pi_{k,d}}\frac{\log\lambda}{16}+\frac{m_k+N_d}{\pi_{k,d}}\big(-\frac{\log{\lambda}}{\alpha}+\frac{\log\lambda}{16}\big)
\\
&>\frac{-3\log\lambda}{4}
\end{align*}

By Pliss lemma,  there exist a number  $\rho\in(0,1)$ only depending on $\lambda$ and a sequence of points   $ \{s_1,\cdots,s_l\}\subset\{0,\cdots, \pi_{k,d}-1\}$ such that
\begin{itemize}
\item $$\frac{l}{\pi_{k,d}}\geq\rho;$$
\item  $$ \prod_{l=j}^{s_i-1}\norm{Df|_{E^c(f^{l-N_d}(z_k))}}
 \geq (\frac{1}{\sqrt{\lambda}})^{s_i-j+1}, \textrm{ for any $j=0,\cdots,s_i-1$}.$$
 \end{itemize}

By   Equation (\ref{equ:average on center}), we have that  $\{s_1,\cdots,s_l\}\cap [N_d+N, N_d+t_k]=\emptyset$. Since the center Lyapunov exponent of $\cO_p$ is $\frac{-\log\lambda}{\alpha}$ and $f^{t_k+N_d}(z_k)\in W^s_{d/2}(\cO_p)$, when   $k$ is chosen large,  we have that
$$\frac{1}{\pi_{k,d}-i}\sum_{j=i}^{\pi_{k,d}-1}\log\norm{Df|_{E^c(f^{j-N_d}(z_k))}}>\frac{-\log\lambda}{2}, \textrm{for any $i=0,\cdots,\pi_{k,d}-1$.}$$

Since $E^{ss}$ is uniformly contracting and $E^{ss}\oplus (E^c\oplus E^{uu})$ is a dominated splitting over $\tilde{\La}$, the orbit segment  $\{f^{-N_d}z_{k},\ldots,f^{N_{d}+t_k+m_{k}}(z_k)\}$ is a $\sqrt{\lambda}$-quasi hyperbolic string corresponding to the splitting
                  $E^{ss}\oplus (E^c\oplus E^{uu})$.

This ends the proof of Claim~\ref{c.choice of t_k}.
\endproof

     By Lemma \ref{shadow}, there exists  a periodic orbit $\cO_{p_{1}}$ of period $\pi_{k,d}$ such that
      $$\ud(f^i(p_1),f^i(z_{k,d}))<L\cdot d, \textrm{ for any $i\in[0, \pi_{k,d}-1]$.}$$

When $d$ is chosen small enough, by the second item of  Claim \ref{c.choice of t_k} and uniform continuity of the function $\log\norm{Df|_{_{E^c}}}$ defined on $\tilde\La$,
      we have that $p_1$ is a $(\frac{1}{\sqrt[4]{\lambda}}, E^c\oplus E^{uu})$ expanding point; By Lemma \ref{lp}, the point $p_1$ has uniform size of unstable manifold independent of $d$.   Once again, when $d$ is chosen small, by the fact that  the strong stable manifolds of $p_1$ and of $\cO_p$ are  the stable manifolds of $p_1$ and of $\cO_p$ respectively, we have that $\cO_{p_1}$ and $\cO_p$ are homoclinically related.

      On the other hand, when $d$ is chosen small, by the first item of Claim \ref{c.choice of t_k} and the first item of
      Lemma ~\ref{l.orbit follow the measure},
      one can check that
$$\ud(\delta_{\cO_{p_1}},\alpha\delta_{\mathcal{O}_{p}}+(1-\alpha)\nu)<4\varepsilon.$$
Hence, $\alpha\delta_{\mathcal{O}_{p}}+(1-\alpha)\nu$ is approximated by
periodic measures whose supports are periodic orbits homoclinically relate to $\cO_p$.

     By the arbitrary choice of $\alpha$ and compactness of the set $\{\alpha\delta_{\mathcal{O}_{p}}+(1-\alpha)\nu|\alpha\in[0,1]\}$, $\nu$ is  approximated by periodic measures, ending the proof of Proposition~\ref{p.main}.
    \endproof

Now we are ready to give the proof of Theorem~\ref{thmF}.
\proof[Proof of Theorem~\ref{thmF}]
Recall that the open set $V$ is a small neighborhood of the split flip flop configuration $(\La, C,\cC^{uu},\mathfrak{D},\cO_q,\Delta^s,\Delta^u)$ such that the maximal invariant set $\tilde\La$  in $\overline{V}$ is partially hyperbolic with center dimension one.  Up to shrinking $V$, we can assume that $V=U \cup V_1\cup V_2\cup W$ such that
\begin{itemize}
\item $U$ is a small  open neighborhood of $C$ satisfying that the maximal invariant set in $\overline{U}$ is $\La$;
\item $V_1$ and $V_2$ are small neighborhoods of $\cup_{i\in\N}f^{i}(\Delta^s)$ and of $\cup_{i\in\N}f^{-i}(\Delta^u)$ respectively;
\item $W$ is a small neighborhood of $\cO_q$ such that the maximal invariant set in $\overline{W}$ is $\cO_q$.
\end{itemize}

Now, we will choose a small neighborhood $V_0\subset V$ of the flip flop configuration such that any non-hyperbolic ergodic measure supported on the maximal invariant set $\tilde\La_0$ in $\overline{V_0}$ satisfying the conditions in Proposition ~\ref{p.main}.

By assumption,
there exists an integer $N$ such that for any point $x\in \cap_{i=-N}^{N}f^{i}(U)\cap\tilde\La$, there exists a periodic point $p\in\La$ such that  $W^{ss}_{loc}(x)$ intersects $W^{u}_{loc}(p)$.
For simplicity, we assume that  the periodic point $q$ is a fixed point. We take a small neighborhood $W^{\prime}\subset W$ of $\cO_q$ such that
$$\log{\norm{Df|_{E^c(x)}}}<\lambda<0, \textrm{ for any point $x\in W^{\prime}\cap \tilde\La$}. $$
 On the other hand, there exists an integer $N_0$ such that $f^{N_0}(\Delta^s)\cup f^{-N_0}(\Delta^u) $ is contained in $W^{\prime}$.
Let $\tilde{N}$ be the smallest integer satisfying $$\tilde{N}>(2N+2N_0)\frac{b}{|\lambda|},\textrm{ where  $b=\max_{x\in M}\norm{Df(x)}$.}$$
 We take a small neighborhood $W_0$ of $\cO_q$ such that
$$W_0\cup\cdots\cup f^{\tilde{N}}(W_0)\subset W^{\prime}.$$
Let $U_0\subsetneq U$ be a neighborhood of $C$.  By the first and the second items in the definition of flip flop configuration, one can take small neighborhoods $V_1^{\prime}\subset V_1$ and $V_2^{\prime}\subset V_2$ of the sets $Orb^{+}(\Delta^s,f)$ and $Orb^{-}(\Delta^u,f)$ respectively such that any point $x\in\tilde\La_0\backslash\La$, where $\tilde\La_0$ is the maximal invariant set of $U_0\cup V_1^{\prime}\cup V_2^{\prime}\cup W_0$, the positive orbit of $x$ intersects $W_0$.

Let $V_0=U_0\cup V_1^{\prime}\cup V_2^{\prime}\cup W_0\subset V$.
By the choice of $W^{\prime}$ and $W_0$, for any  non-hyperbolic ergodic measure $\nu\in\cM_{erg}(\tilde\La_0,f)$ and any point $x$ in the basin of $\nu$, the forward orbit of $x$ contains an orbit segment of length $2N+1$ which are contained in $U$, hence there exists a point $y$ in the basin of $\nu$ such that $W^{ss}_{loc}(y)$  intersects the local unstable manifold of a periodic point contained in $\La$ in a point $z$; moreover, by the uniform contraction for the local strong stable manifold,    the closure of the orbit of $z$ is strictly contained in $V$.  Now, by applying $\nu$, $\tilde\La$ and $V$ to Proposition~\ref{p.main}, one has that $\nu$ is accumulated by periodic measures, ending the proof of Theorem~\ref{thmF}.
\endproof

\section{The closure of periodic measures contains a segment joining $\delta_{\cO_q}$ to a measure in the blender: Proof of Theorem~\ref{thmE}}
Given a split flip flop configuration $(\Lambda, C, \mathcal{C}^{uu},\mathfrak{D},\cO_q, \Delta^s,\Delta^u)$ of $f\in\diff^1(M)$ formed by a  blender horseshoe $(\Lambda, C, \mathcal{C}^{uu},\mathfrak{D})$ and a hyperbolic periodic point $q$. In this section, we prove that  there exists an invariant measure $\mu$ (maybe non-ergodic) supported on $\Lambda$ such that the convex combination $\{\alpha\mu+(1-\alpha)\delta_{\cO_q};\alpha\in[0,1]\}$ is approximated  by periodic measures.

We take a small neighborhood $V$ of the split flip flop configuration such that the maximal invariant $\tilde\La$ in $\overline{V}$ is partially hyperbolic with center dimension one. We assume, in addition, that there exist two $Df$ strictly  invariant  cone fields $\cC^u_V$ and $\cC^{uu}_V$ in $V$, which are  continuous extensions of the center unstable cone fields $\cC^u$ and the strong unstable cone field $\cC^{uu}$ in $C$ respectively.

   Let's fix a sequence of function $\{g_{i}\}_{i=1}^{+\infty}$ which is a dense subset of $C^{0}(M,\R)$. Then $\{g_{i}\}_{i=1}^{+\infty}$
  determines a metric on the probability measure space on $M$ in the following way: for any probability measures $\nu_1,\nu_2$ on $M$,
  we have
  \begin{displaymath}
 \ud(\nu_1,\,\nu_2)=\sum_{i=1}^{\infty}\frac{|\int g_i\,\ud\nu_1-\int g_i\,\ud\nu_2|}{2^{i}\norm{g_i}_{_{C^0}}}.
\end{displaymath}

Since the disc $\Delta^u\subset W^u(q)$  belongs to $\mathfrak{D}$, by the strictly invariant property of $\mathfrak{D}$, this segment would intersect with the local stable manifold of $\Lambda$ in a (Cantor) set which is denoted as $\mathfrak{C}$.
For any point $x\in\mathfrak{C}$,  there exists a sequence of discs $\{D_i\}_{i\in\N}\subset\mathfrak{D}$
such that
\begin{itemize}
\item $f^i(x)\in D_i$ for any $i\in\N$;
\item $D_0=\Delta^u$ and $D_{i+1}\subset f(D_i)$, for any $i\in\N$.
\end{itemize}

\begin{theorem} \label{convex}With the assumption above. Given  $x\in\mathfrak{C}$ and let $\mu$ be an  accumulation of $\{\frac{1}{n}\sum_{i=0}^{n-1}\delta_{f^{i}(x)}\}_{n\in\N}$.
 Then the convex combination $\{\alpha\mu+(1-\alpha)\delta_{\mathcal{O}_{q}}|\,\alpha\in[0,1]\}$ is contained in
 the closure of the set of periodic measures.
\end{theorem}

Now Theorem~\ref{thmE} is directly from Theorem~\ref{convex}. Hence, we only need to prove Theorem~\ref{convex}.

Recall that 
 the mean center Lyapunov exponent of the invariant measure $\alpha\mu+(1-\alpha)\delta_{\mathcal{O}_{q}}$, by definition, is
 $$\lambda^{c}(\alpha\mu+(1-\alpha)\delta_{\mathcal{O}_{q}})=\alpha \lambda^{c}(\mu)+(1-\alpha)\lambda^{c}(\delta_{\mathcal{O}_{q}}),$$
   hence there exists $\alpha_{0}\in[0,1]$ such that
    $$\alpha_{0} \lambda^{c}(\mu)+(1-\alpha_{0})\lambda^{c}(\delta_{\mathcal{O}_{q}})=0.$$

    The proof of  Theorem~\ref{convex} consists in two parts.  We first show that for any $\alpha\in[0,\alpha_0]$, $\alpha\mu+(1-\alpha)\delta_{\mathcal{O}_{q}}$ is accumulated by periodic measures.
    Then we show the other half convex combination is also approached by periodic measures.
    The proof of these two parts are quite different. But  their proofs still consist in  finding  quasi hyperbolic periodic pseudo orbits and applying Lemma \ref{shadow} to find the periodic orbits.

  \begin{lemma} \label{halfconvex}
  For any $\alpha\in[0,\alpha_{0}]$, the invariant measure $\alpha\mu+(1-\alpha)\delta_{\mathcal{O}_{q}}$ is  accumulated by a sequence of periodic orbits which are homoclinically related to $\mathcal{O}_{q}$ in $V$.
  \end{lemma}
  \proof We fix $\alpha\in(0,\alpha_{0})$, then  $\lambda^{c}(\alpha\mu+(1-\alpha)\delta_{\mathcal{O}_{q}})$ is negative.

  Let $\lambda=\exp(\lambda^{c}(\alpha\mu+(1-\alpha)\delta_{\mathcal{O}_{q}}))$.
  By Lemma \ref{shadow}, there exist two positive numbers $L$ and $d_{0}$ such that for any $d\in(0,d_0]$, every $(\lambda+1)/2$-quasi hyperbolic periodic $d$-pseudo orbit  corresponding  to the splitting $T_{\tilde\La}M=(E^{ss}\oplus E^c)\oplus E^{uu}$ is $L\cdot d$ shadowed by a periodic orbit.

   By the continuity of center distribution, there exists $\delta>0$ such that for any $z_1, z_2\in\tilde{\Lambda}$ satisfying $\ud(z_1, z_2)<\delta$,
    we have that
  $$\frac{4\lambda}{1+3\lambda}\leq\frac{\norm{Df|_{E^{c}(z_1)}}}{\norm{Df|_{E^{c}(z_2)}}}\leq \frac{1+3\lambda}{4\lambda}.$$
For any  $d\in(0,\min\{d_{0},\frac{\delta}{L}\})$  whose precise value would be fixed at the end,   there exists a positive integer $N_{d}$ such that
 $$f^{N_{d}}(\Delta^{s})\subset W^s_{d/2}( q) \textrm{ and }  f^{-N_d}(\Delta^{u})\subset W^u_{d/2}(q).$$

Let $\tau_0>1$ be a number such that for any point $x\in V\cap f^{-1}(V)$ and $v\in\cC^{uu}_V(x)$, one has $\norm{Df(v)}\geq \tau_0\norm{v}$.
We denote by $N_{\delta}=[\frac{\log{b_0}-\log\delta}{\log\tau_0}]+1$,
where $b_0$ is an upper bound for the diameters of the discs in $\mathfrak{D}$.

   In the following, we will find a $(\lambda+1)/2$-quasi hyperbolic periodic $d$-pseudo orbit which will  stay almost $\alpha$ proportion of time to follow an orbit segment of $x$ and $(1-\alpha)$ proportion of time to follow the orbit of $\mathcal{O}_{q}$; then we apply Lemma \ref{shadow}.

 For any $\epsilon>0$, there exists an  integer  $n$ arbitrarily  large such that
 $$\ud(\frac{1}{n}\sum_{i=0}^{n-1}\delta_{f^{i}(x)},\mu)<\epsilon.$$
   We choose the $\delta$ neighborhood of $x$ in $W^{u}(\mathcal{O}_{q})$ and denote it as $D_x^{u}$.
  Consider the connected component $D^u_x(n)$ of $f^{n}(D_x^{u})\cap B_{\delta}(f^{n}(x))$ which contains $f^{n}(x)$.
   By the choice of $N_{\delta}$, one gets that  $f^{N_{\delta}}(D^u_x(n))$ contains a disc in $\mathfrak{D}$, hence  $f^{N_{\delta}}(D^u_x(n))$ has transverse intersection with  $\Delta^{s}$. By the choice of $N_{d}$, there exists a transverse intersection $y$ between $f^{N_{\delta}+N_{d}}(D^u_x(n))$ and  $W^{s}_{d/2}(q)$.

    Consider the orbit segment  $\sigma_{m,n}=\{f^{-m\pi(q)-n-2N_{d}-N_{\delta}}(y), \cdots, y\}$. We denote by
    $$t_{m,n}=m\pi(q)+n+2N_{d}+N_{\delta} \textrm{ and } t_n= n+2N_{d}+N_{\delta}.$$
    Notice that for any $m\in\N$, we have $\ud(f^{t_{m,n}}(y),y)<d.$

    We denote by $b=\max\{\sup_{x\in M}\norm{Df(x)},\sup_{x\in M}\norm{Df^{-1}(x)}\}$.
    \begin{claim}\label{c.choice of n and m} There exist $n$ and $m$ arbitrarily large such that
    \begin{itemize}
     \item $$\Big|\frac{n}{m\pi(q)}-\frac{\alpha}{1-\alpha}\Big|+\frac{2N_{d}+N_{\delta}}{n}\cdot b<\epsilon;$$
    \item $$\ud(\frac{1}{t_n}\sum_{i=0}^{t_n-1}\delta_{f^{-i}(y)},\mu)<\epsilon;$$
    \item   $\sigma_{m,n}$ is a $(1+\lambda)/2$-quasi hyperbolic string corresponding to the splitting $(E^{ss}\oplus E^c)\oplus E^{uu}$.
    \end{itemize}
    \end{claim}
    The proof of Claim \ref{c.choice of n and m} is  just  like the one of Claim \ref{c.choice of t_k}.

    Using Lemma \ref{shadow}, we can get a periodic orbit $\cO_{q_1}$ of period $t_{m,n}$ such that
     $$\ud(f^j(q_1),f^j(f^{-t_{m,n}}(y)))<L\cdot d, \textrm{for any $j=0,\cdots, t_{m,n}-1$}.$$
      Arguing as before, when $d$ is chosen small enough, one has  that
       \begin{itemize}
       \item  $\ud(\delta_{\cO_{q_1}},\alpha\mu+(1-\alpha)\delta_{\mathcal{O}_{q}})<4\cdot\epsilon;$
         \item $\cO_{q_1}$ is homoclinically related to $\cO_q$.
        \end{itemize}
        Hence, $\alpha\mu+(1-\alpha)\delta_{\mathcal{O}_{q}}$ is approximated by periodic orbits which are homoclinically related to $\mathcal{O}_{q}$.

    This ends the proof of Lemma \ref{halfconvex}.
 \endproof

   For $\alpha\in[0, \alpha_0)$,  the property of  $\Delta^{u}$ helps to find the type of quasi hyperbolic string that we need.
   For the case $\alpha\in(\alpha_0,1]$,  the quasi hyperbolic string that we need  is another type, that is,  we want $Df$ along the  center direction to have expanding behavior on the quasi hyperbolic string.
Indeed, using the strategy above, we can start  from a small neighborhood of $\La$ then go arbitrarily close to $\cO_q$ to stay  for arbitrarily long time; however,  after that, it is  not clear  if we can go    arbitrarily close to  $x$ by an arbitrarily  small proportion of time.

    To deal with this situation, we change the strategy. The proof for the case $\alpha\in(\alpha_0, 1]$  strongly depends on the fact that
    $\alpha_0\mu+(1-\alpha_0)\delta_{\mathcal{O}_{q}}$ is approximated by hyperbolic periodic orbits homoclinically related to $\cO_q$.


   \begin{proposition} \label{gap} With the assumption we posed at the beginning of this section. There exists a constant $\rho>0$, such that for any hyperbolic periodic orbit $\cO_{q^{\prime}}$ which is homoclinically related to $\cO_q$ in $V$,
   any $\epsilon>0$, any hyperbolic  periodic orbit $\mathcal{O}_{p}\subset\Lambda$ and $\alpha\in(0,1)$ satisfying that $$\lambda^{c}\big(\alpha\delta_{\mathcal{O}_{p}}+(1-\alpha)\,\delta_{\mathcal{O}_{q^{\prime}}}\big)>0,$$
one has that there exists a hyperbolic periodic orbit $\mathcal{O}_{p_0}$ homoclinically related to $\cO_p$ such that: $$\ud\big(\delta_{\mathcal{O}_{p_0}},\,\alpha\delta_{\mathcal{O}_{p}}+(1-\alpha)\,\delta_{\mathcal{O}_{q^{\prime}}}\big)
   <\rho\cdot (1-\alpha)\cdot|\lambda^c({\cO_{q^{\prime}}})|+\epsilon.$$
   \end{proposition}
   \proof  For $\epsilon>0$, there exists  an integer  $N$ large enough such that $$2\sum_{i=N}^{\infty}\frac{1}{2^{i}}<\frac{\epsilon}{2}.$$
   Then there exists $\delta>0$ such that for any $x,y$ satisfying $\ud(x,y)<\delta$, we have that $$|g_i(x)-g_i(y)|<\frac{\epsilon}{8} \norm{g_i}_{_{C^0}}, \textrm{ for any $i=1,\cdots,N$}.$$

  Let $\epsilon_0$ be the strength of $\mathfrak{D}$.
  Since $\mathcal{O}_{q^{\prime}}$ is homoclinically related to  $\mathcal{O}_q$ in $V$,  by Inclination Lemma and the definition of flip-flop configuration,   there exist  two compact submanifolds $\Delta^s(q^{\prime})\subset W^{s}(\mathcal{O}_{q^{\prime}})$ and $\Delta^u(q^{\prime})\subset W^{u}(\mathcal{O}_{q^{\prime}})$ such that
  \begin{itemize}
    \item $$\Delta^u(q^{\prime})\in\cV_{\epsilon_0/4}(\mathfrak{D});$$
    \item    $$Orb^{-}(\Delta^{u}(q^{\prime}))\cup Orb^{+}(\Delta^{s}(q^{\prime}))\subset V;$$

  \item  for any  $D^{u}\in\mathfrak{D}$,  the disc $D^{u}$
  intersects the interior of  $\Delta^{s}(q^{\prime})$  transversely.
  \end{itemize}

    Let $\lambda=\exp\big(-\lambda^{c}(\alpha\delta_{\mathcal{O}_{p}}+(1-\alpha)\delta_{\mathcal{O}_{q^{\prime}}})\big)$. By Lemma \ref{shadow}, there exist two positive numbers $L$ and $d_{0}$ such that for any $d\in(0,d_0]$, every $(\lambda+1)/2$-quasi hyperbolic periodic $d$-pseudo orbit   corresponding to the splitting $T_{\tilde\La}M=E^{ss}\oplus (E^c\oplus E^{uu})$ is $L\cdot d$ shadowed by periodic orbit.

    We choose a number $d\in(0,\min\{d_{0},\delta\}]$  small enough such that $L\cdot d$ is much less than $\delta$; the precise value of $d$ would be fixed at the end.   By the choice of $\Delta^{s}({q^{\prime}})$ and $\Delta^{u}({q^{\prime}})$,  there exists an integer $N^1_{d}$ satisfying that:
 $$f^{N_{d}^1}(\Delta^{s}(q^{\prime}))\cup f^{-N_{d}^1}(\Delta^{u}(q^{\prime}))\subset B_{d/2}(\cO_{q^{\prime}}).$$
 Up to increasing $N_d^1$, we can assume that $f^{N_{d}^1}(W^{s}_{loc}(\Lambda))\subset W^{s}_{d/2}(\Lambda)$.

By the transitivity of $\La$,  Remark \ref{r.domain of blender} and the strictly invariant property of $\mathfrak{D}$, there exist  $N^2_d\in\N$ and a disc   $D_0\subset W^{u}_{d/2}(\cO_p)$ such that
  $$f^{N^2_{d}}(D_0)\in\mathfrak{D} \textrm{ and } \cup_{i=0}^{N^2_d}f^i(D_0)\subset C.$$
   Up to increasing $N^1_d$ or $N^2_d$, we can assume that $N_d^1=N_d^2$ and we denote it by $N_d$.

 We fix two  plaque families $W^{cs}$ and $W^{cu}$ corresponding to the bundle $E^{ss}\oplus E^c$ and to  bundle $E^c\oplus E^{uu}$ respectively.

By the choice of $\Delta^s(q^{\prime})$,  the disc $f^{N_d}(D_0)$ intersects $\Delta^{s}(q^{\prime})$ transversely and we denote the intersection as $y$, then one has  that
  \begin{itemize}\item $$\overline{Orb(y,f)}\subset V;$$
  \item   $$Orb^{-}(f^{-N_d}(y))\subset W^{u}_{d/2}(\cO_p) \textrm{  and } Orb^{+}(f^{N_d}(y))\subset W^{s}_{d/2}(\cO_{q^{\prime}}).$$
    \end{itemize}
     By Lemma \ref{lp}, there exists a number $\delta_d>0$ such that
     $$W^{cs}_{2\delta_d}(y)\subset f^{-N_d}(W^s_{d/2}(\cO_q^{\prime}))
      \textrm{ and } W^{cu}_{2\delta_d}(y)\subset f^{N_d}(W^{u}_{d/2}(\cO_p)).$$

Since $\Delta^{u}(\cO_{q^{\prime}})\subset W^{u}(\cO_{q^{\prime}})$,
by uniform expansion of $Df$ along the strong unstable cone field $\cC^{uu}_V$,
there exists an integer $N_d^{\prime}$ large such that
for any disc $D$ tangent to the  cone field $\cC^{uu}_V$,
 if $D$ intersects $W^s_{d/2}(\cO_{q^{\prime}})$ transversely in a point whose distance to the relative boundary of $D$ is no less than $\delta_d$,
 we have that $f^{N_d^{\prime}}(D)$ contains a disc belong to $\cV_{\epsilon_0/2}(\mathfrak{D})$.
  Up to increasing $N_d$ or $N_d^{\prime}$, we can assume that $N_d=N_d^{\prime}$.

Let $S(y)$ be the $\delta_d$ tubular neighborhood of $W^{uu}_{\delta}(y)$ in $W^{cu}(y)$, hence one has  $S(y)\subset W^{cu}_{2\delta_d}(y)$.
We denote by  $S_n(y)$ the connected component  of $f^{n+N_d}(S(y))\cap B_{d/2}(\cO_q^{\prime})$ which contains $f^{n+N_d}(y)$ for any  $n\in\N$, and we denote by \\ $b=\max\{\sup_{x\in M}\norm{Df_x},\sup_{x\in M}\norm{Df^{-1}_x}\}$.
Since the center Lyapunov exponent of  the orbit of $q^{\prime}$ is negative, when $n$ is chosen large enough, arguing as Claim~\ref{c.length of central curve} and Claim~\ref{c.lower bound of central length},
one has that  $S_n(y)$  is  a uu-foliated cu-disc  satisfying that
 \begin{itemize}
 \item the  central length of $S_n(y)$   is at least $\exp(n(\lambda^c(\cO_{q^{\prime}})-\epsilon))\cdot \delta_d \cdot b^{-N_d}$;
  \item $S_n(y)$ is foliated by discs of size $\delta_d$ tangent to the cone field $\cC^{uu}_V$.
 \end{itemize}
 Then $f^{N_d}(S_n(y))$ contains a cu-strip $D^{cu}$ in $C$ foliated by uu-discs  in $\mathcal{V}_{\epsilon_0/2}(\mathfrak{D})$ such that
  for the central length $\ell^c(D^{cu})$ of $D^{cu}$, one has
 $$ \ell^c(D^{cu})\geq \exp(n\lambda^c(\cO_{q^{\prime}}))\cdot \delta_d\cdot b^{-N_d}.$$
  Let $\tau>1$ be a number such that
  $$\norm{Df(v)}\geq\tau\norm{v}, \textrm{ for any $x\in C\cap f^{-1}(C)$ and any $v\in\cC^u(x)$.}$$
  Following the strategy   in Lemma \ref{l.orbit follow the measure},
   there exists a constant $c_d$ independent of $n$ and an integer $k$ such that
    \begin{itemize}
   \item $f^{k}(D^{cu})$ intersects the local stable manifold of $\mathcal{O}_p$;
   \item $f^{i}(D^{cu})$ does not intersect the local stable manifold of $\mathcal{O}_p$, for $i<k$;
\item  when $n$ is chosen large enough, we have an upper bound for $k$:
 \begin{equation}\label{equ:upper bound of k}  k\leq \frac{2n\cdot\pi(q^{\prime})\cdot |\lambda^{c}(\mathcal{O}_{q^{\prime}})|+c_d}{\log\tau}.
 \end{equation}
 \end{itemize}
 Let $x$ be the intersection of $f^{k}(D^{cu})$ and $W^{s}_{loc}(\cO_p)$.

  We denote by $x_{n,d}=f^{-2N_{d}-n\pi(q^{\prime})-k}(x).$  For any positive integer $m$, consider the orbit segment
   $$\sigma_{m,n}=\{x_{n,d},\cdots,f^{m\pi(p)+N_{d}}(x)\},$$
   Notice that $\ud(x_{n,d}, f^{m\pi(p)+N_{d}}(x))<d.$
   \begin{claim}\label{c.gap}There exist  integers $m$ and $n$ arbitrarily large such that
   \begin{itemize}
     \item $$\big|\frac{n\pi(q^{\prime})}{m\pi(p)}-\frac{1-\alpha}{\alpha}\big|
         +\big|\frac{m\pi(p)}{n\pi(q^{\prime})}-\frac{\alpha}{1-\alpha}\big|<\frac{\epsilon}{16};$$
         \item $$\frac{c_d+3N_d}{n}\cdot (1+\frac{1}{\log\tau})<\frac{\epsilon}{16}$$
     \item $\{x_{n,d},\cdots,f^{m\pi(p)+N_{d}}(x)\}$ is a $(1+\lambda)/2$ quasi hyperbolic string corresponding to the splitting $E^{ss}\oplus (E^c\oplus E^{uu})$.
    \end{itemize}
    \end{claim}
      \begin{figure}[h]
\begin{center}
\def\svgwidth{0.9\columnwidth}
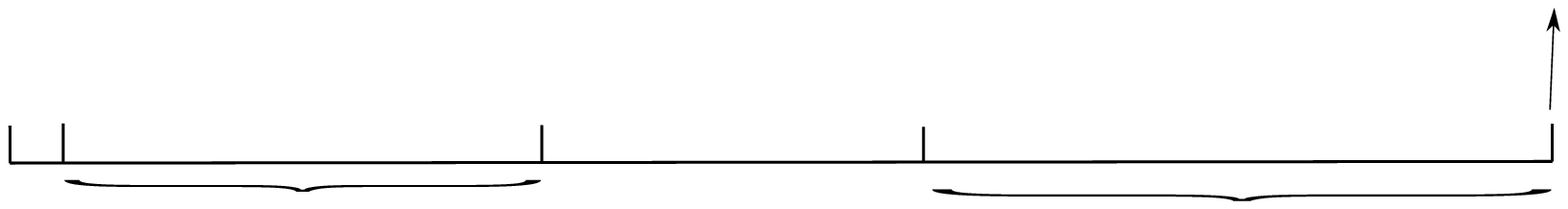
  \caption{}
\end{center}
\end{figure}

    The proof of this claim is exactly as the proof of Claim~\ref{c.choice of  t_k}.

    By Lemma \ref{shadow}, there exists a periodic orbit $\mathcal{O}_{p_{_0}}$  of period $n\pi(q^{\prime})+m\pi(p)+3N_{d}+k$ such that for any $i=0, \cdots, \pi(p_{_0})-1 $, we have
    $$\ud(f^{i}(x_{n,d}),\,f^{i}(p_{_0}))<L\cdot d.$$
    Arguing as before, when $d$ is chosen small enough, we have that $\mathcal{O}_{p_{_0}}$ is homoclinically related  to $\cO_p$ in $V$.

     We denote by
    \begin{itemize}\item[--] $I_1=\{N_d, \cdots, N_d+n\pi(q^{\prime})-1\}$;
    \item[--] $I_2=\{0,\cdots, N_d-1\}\cup\{N_d+n\pi(q^{\prime}), \cdots, 3N_d+k+n\pi(q^{\prime})-1\}$;
     \item[--]$I_3=\{3N_{d}+k+n\pi(q^{\prime}), \cdots,\pi(p_0)-1\}$.
    \end{itemize}
    By the choice of $\sigma_{m,n}$ and the fact that  the orbit of $p_0$ shadows $\sigma_{m,n}$ at a distance of $L\cdot d$,  we have that
    \begin{itemize}
    \item
    $$\ud(f^j(p_0), f^{j-N_d}(q^{\prime}))<L\cdot d+d, \textrm{ for any $j\in I_1$};$$
    \item
    $$\ud(f^j(p_0), f^{j-3N_{d}-k-n\pi(q^{\prime})}(p))<L\cdot d+d, \textrm{ for  any $j\in I_3$}.$$
    \end{itemize}

    \begin{claim}\label{c.integration distance}For each integer $i\in[1,N]$, we have that
     $$\Big|\int g_i\,\ud\delta_{\mathcal{O}_{p_{_0}}}-\alpha\int g_i\,\ud\delta_{\mathcal{O}_{p}}-(1-\alpha)\int g_i\,\ud\delta_{\mathcal{O}_{q^{\prime}}}\Big|
    <\big(\frac{\epsilon}{2}+\frac{2|\lambda^{c}(\cO_{q^{\prime}})|}{\log\tau}\cdot(1-\alpha)\big)\cdot \norm{g_i}_{_{C^0}}.$$
    \end{claim}
    \proof For each $i\in[1,N]$, we have that
    \begin{align*}
    &\Big|\int g_i\,\ud\delta_{\mathcal{O}_{p_{_0}}}-\alpha\int g_i\,\ud\delta_{\mathcal{O}_{p}}-(1-\alpha)\int g_i\,\ud\delta_{\mathcal{O}_{q^{\prime}}}\Big|
    \\
    &\leq \Big|\frac{1}{\pi(p_{_0})}\sum_{j\in I_1}g_{i}(f^{j}(p_{_0}))-(1-\alpha)\int g_i\,\ud\delta_{\mathcal{O}_{q^{\prime}}}\Big|+\Big|\frac{1}{\pi(p_{_0})}\sum_{j\in I_2}g_{i}(f^{j}(p_{_0}))\Big|
    \\
    &\hspace{5mm}+
   \Big|\frac{1}{\pi(p_{_0})}\sum_{j\in I_3}g_{i}(f^{j}(p_{_0}))-\alpha \int g_i\,\ud\delta_{\mathcal{O}_{p}}\Big|.
    \end{align*}

    By the choice of $I_1$,   we have the estimate:
    \begin{align*}
    &\Big|\frac{1}{\pi(p_{_0})}\sum_{j\in I_1}g_{i}(f^{j}(p_{_0}))-(1-\alpha)\int g_i\,\ud\delta_{\mathcal{O}_{q^{\prime}}}\Big|
    \\
    &<\Big|\frac{1}{\pi(p_{_0})}\big(\sum_{j\in I_1}g_{i}(f^{j}(p_{_0}))-\sum_{j\in I_1}g_{i}(f^{j}(q^{\prime})\big)\Big|
    +\Big|(\frac{n\pi(q^{\prime})}{\pi(p_{_0})}-1+\alpha)\int g_i\,\ud\delta_{\mathcal{O}_{q^{\prime}}}\Big|
    \\
    &\leq\frac{\epsilon}{8}\norm{g_i}_{_{C^0}}+\Big|\frac{n\pi(q^{\prime})}{\pi(p_{_0})}-1+\alpha\Big|\cdot\norm{g_i}_{_{C^0}}
    \end{align*}
    By Equation (~\ref{equ:upper bound of k}) and Claim ~\ref{c.gap}, when $m,n$ are chosen large,  we have that
    \begin{align*}
    \big|\frac{n\pi(q^{\prime})}{\pi(p_{_0})}-1+\alpha\big|
    &\leq\big|\frac{n\pi(q^{\prime})}{\pi(p_{_0})}-\frac{n\pi(q^{\prime})}{n\pi(q^{\prime})+m\pi(p)}\big|
    +\big|\frac{n\pi(q^{\prime})}{n\pi(q^{\prime})+m\pi(p)}-1+\alpha\big|
    \\
    &\leq  \frac{(k+3N_d)}{n\pi(q^{\prime})+m\pi(p)}\cdot \frac{n\pi(q^{\prime})}{n\pi(q^{\prime})+m\pi(p)}+\frac{\epsilon}{16}
    \\
    &<\frac{3|\lambda^{c}(\cO_{q^{\prime}})|}{\log\tau}\cdot (1-\alpha)^2+\frac{\epsilon}{8}.
    \end{align*}

    Hence, we have that
    $$\Big|\frac{1}{\pi(p_{_0})}\sum_{j\in I_1}g_{i}(f^{j}(p_{_0}))-(1-\alpha)\int g_i\,\ud\delta_{\mathcal{O}_{q^{\prime}}}\Big|<\frac{\epsilon}{4}\norm{g_i}_{_{C^0}}+\frac{3|\lambda^{c}(\cO_{q^{\prime}})|}{\log\tau}\cdot (1-\alpha)^2\cdot \norm{g_i}_{_{C^0}}$$
    Similarly, by choosing $m,n$ large enough, we also have that
    \begin{align*}&\Big|\frac{1}{\pi(p_{_0})}\sum_{j\in I_2}g_{i}(f^{j}(p_{_0}))\Big|<\frac{3N_d+k}{\pi(p_{_0})}\norm{g_i}_{_{C^0}}
    <\frac{3|\lambda^{c}(\cO_{q^{\prime}})|}{\log\tau}\cdot(1-\alpha)\norm{g_i}_{_{C^0}},
      \\
            &\Big|\frac{1}{\pi(p_{_0})}\sum_{j\in I_3}g_{i}(f^{j}(p_{_0}))-\alpha\int g_i\,\ud\delta_{\mathcal{O}_{p}}\Big|
              <\frac{\epsilon}{4}\norm{g_i}_{_{C^0}}+\frac{3|\lambda^{c}(\cO_{q^{\prime}})|}{\log\tau}\cdot \alpha\cdot(1-\alpha)\cdot\norm{g_i}_{_{C^0}}
    \end{align*}
    Hence, we have that for any $i\in[1,N]$,
    $$\Big|\int g_i\,\ud\delta_{\mathcal{O}_{p_{_0}}}-\alpha\int g_i\,\ud\delta_{\mathcal{O}_{p}}-(1-\alpha)\int g_i\,\ud\delta_{\mathcal{O}_{q^{\prime}}}\Big|
    <\big(\frac{\epsilon}{2}+\frac{6|\lambda^{c}(\cO_{q^{\prime}})|}{\log\tau}\cdot(1-\alpha)\big)\cdot \norm{g_i}_{_{C^0}}.$$

   This ends the proof of  Claim~\ref{c.integration distance}
   \endproof

    By the choice of $N$ and Claim \ref{c.integration distance}, we have
    \begin{align*} \ud\big(\delta_{\mathcal{O}_{p_{_{_0}}}},\,\alpha\delta_{\mathcal{O}_{p}}+(1-\alpha)\delta_{\mathcal{O}_{q^{\prime}}}\big)
    &=\sum_{i=1}^{\infty}\frac{|\int g_i\,\ud\delta_{\mathcal{O}_{p_{_0}}}-\alpha\int g_i\,\ud\delta_{\mathcal{O}_{p}}-(1-\alpha)\int g_i\,\ud\delta_{\mathcal{O}_{q^{\prime}}}|}{2^{i} \norm{g_i}_{_{C^0}}}
    \\
    &<\sum_{i=1}^{N}\frac{|\int g_i\,\ud\delta_{\mathcal{O}_{p_{_0}}}-\alpha\int g_i\,\ud\delta_{\mathcal{O}_{p}}-(1-\alpha)\int g_i\,\ud\delta_{\mathcal{O}_{q^{\prime}}}|}{2^{i} \norm{g_i}_{_{C^0}}}+\frac{\epsilon}{2}
    \\
    &<\epsilon+\frac{6|\lambda^{c}(\cO_{q^{\prime}})|}{\log\tau}\cdot(1-\alpha).
    \end{align*}

    We take $\rho=\frac{6}{\log\tau}$, ending the proof of Proposition \ref{gap}.
   \endproof
   \begin{Remark} The conclusion of Proposition \ref{gap} also explains the main obstruction to obtain the approximation of  the convex combination between two hyperbolic ergodic measures of different indices by periodic measures.
   \end{Remark}

As an application of Proposition \ref{gap}, we have the following corollary:
\begin{corollary}\label{gap convex} For any $\alpha\in(\alpha_{0},1]$,  the measure $\alpha\mu+(1-\alpha)\delta_{\mathcal{O}_{q}}$ is accumulated by a sequence of periodic orbits which are homoclinically related to $\Lambda$.
   \end{corollary}
    \proof
    By Lemma \ref{halfconvex}, we have a sequence of hyperbolic periodic orbits $\mathcal{O}_{q_{n}}$ which are related to $\mathcal{O}_{q}$ such that $\delta_{\mathcal{O}_{q_{n}}}$
converges to $\alpha_{0}\mu+(1-\alpha_{0})\delta_{\mathcal{O}_{q}}$.

By Theorem \ref{sigmund}, there exists a sequences of periodic orbits $\cO_{p_n}\subset \La$ such that $\delta_{\cO_{p_n}}$ converges to $\mu$.

We denote by
$$B_n=\overline{\{\beta\in[0,1]: \lambda^c(\beta\delta_{\cO_{q_n}}+(1-\beta)\delta_{\cO_{p_n}})>0\}}.$$
By Proposition \ref{gap}, there exists a constant $\rho>0$ such that for any $\beta\in B_n$,
we have that $\beta\delta_{\cO_{q_n}}+(1-\beta)\delta_{\cO_{p_n}}$ is approximated by periodic measures with an error bounded by
 $\rho\cdot \lambda^{c}(\mathcal{O}_{q_n})$.

Since the set $\{\beta\delta_{\cO_{q_n}}+(1-\beta)\delta_{\cO_{p_n}}: \beta\in B_n\}$
tends to the set $\{\alpha\mu+(1-\alpha)\delta_{\mathcal{O}_{q}}: \alpha\in[\alpha_0, 1]\}$
and  $\lambda^c(\cO_{q_n})$ tends to zero,
the invariant measure $\alpha\mu+(1-\alpha)\delta_{\mathcal{O}_{q}}$ is approximated by periodic measures, for any $\alpha\in[\alpha_0,1]$.
\endproof

Now,   Theorem~\ref{convex} is directly from Lemma \ref{halfconvex} and Corollary \ref{gap convex}.
\section{Non-hyperbolic ergodic measures approaching a non-ergodic measure  with vanishing mean center Lyapunov exponent}
We fix a sequence of continuous functions  $\{g_{_i}\}_{i\in\N}\subset C^{0}(M,\R)$, which
  determines a metric $\ud(\cdot,\cdot)$ on the probability measure space  on $M$: for any probability measure $\mu$ and $\nu$, we have the distance:
  $$\ud(\mu,\nu)=\sum_{i=1}^{\infty}\frac{|\int g_{_i}\,\ud\mu-\int g_{_i}\,\ud\nu|}{2^{i}\norm{g_{_i}}_{_{C^0}}}.$$

   Let $f\in\diff^1(M)$, consider a split flip flop configuration  formed by a dynamically defined cu-blender $(\Lambda,U,\cC^{uu},\mathfrak{D})$ and a hyperbolic periodic orbit $\cO_q$. We fix a small  neighborhood $V$ of the split flip flop configuration such that the maximal invariant set $\tilde\La$  of $\overline{V}$ is partially hyperbolic with center dimension one.
  The following result allows us to get a sequence of periodic orbits in $V$ satisfying the \cite{GIKN} criterion.
 \begin{lemma}\cite[Lemma 4.1]{BZ}\label{descend} With the notation above. There exist two constants $\rho>0$ and $\zeta\in(0,1)$, such that  for any $\epsilon>0$ and  any hyperbolic periodic orbit $\gamma$  which   is homoclinically related to $\mathcal{O}_{q}$   inside $V$,
    there exists a hyperbolic  periodic orbit $\gamma^{\prime}$ which is homoclinically related to $\gamma$ in $V$
   satisfying:
  \begin{itemize}
   \item $\lambda^{c}(\gamma^{\prime})>\zeta\lambda^c(\gamma)$;
  \item $\gamma^{\prime}$ is a $(\epsilon, 1-\rho\cdot|\lambda^c(\gamma)|)$ good approximation for $\gamma$.
  \end{itemize}
  \end{lemma}
  \begin{Remark}
  \begin{enumerate}
  \item This result is true for any small neighborhood of the split flip flop configuration;
  \item If $f$ is   partially hyperbolic with center dimension one, one can replace  $V$ by $M$.
  \end{enumerate}
  \end{Remark}
In this section, we show that any invariant measure supported on $\tilde{\La}$, which is  approached by hyperbolic periodic measures of certain index and exhibits vanishing mean center Lyapunov exponent,
is approached by non-hyperbolic ergodic measures.
To be precise, we prove the following:

  \begin{proposition}\label{aa} With the notation above.  Given $\mu\in \mathcal{M}_{inv}(\tilde{\Lambda},f)$  such that  $$\int \log\norm{Df|_{E^{c}}}\ud\mu=0.$$ Assume that  $\mu$ is accumulated by periodic measures whose support are periodic orbits homoclinically related to $\mathcal{O}_q$ inside $V$.

  Then $\mu$ is approximated by non-hyperbolic ergodic measures.
  \end{proposition}
  \proof
  Let $\{\mathcal{O}_{p_{_n}}\}_{n\in\N}$ be a sequence of periodic orbits such  that
  $\delta_{\mathcal{O}_{p_{_n}}}$ converge to $\mu$ and
  $\cO_{p_{_n}}$ is homoclinically related to $\cO_q$ inside $V$.
  We denote by $\lambda_{n}$ the center Lyapunov exponent of $\mathcal{O}_{p_{n}}$,
   then $\lambda_n$ tends to $0$.

  Using Lemma \ref{descend}   and \cite{GIKN} criterion, we will prove that there exists a constant $c>0$, such that  for each  periodic measure $\delta_{\cO_{p_{_n}}}$,  there exists a non-hyperbolic ergodic measure $\nu_n$ satisfying:
  $$\ud(\nu_n,\delta_{\cO_{p_{_n}}})<c\cdot|\lambda_n|.$$

   We fix the periodic orbit $\mathcal{O}_{p_{n}}$, then there exists   an integer $N$ large such that
 $$4\sum_{i=N+1}^{\infty}\frac{1}{2^{i}}\leq |\lambda_{n}|.$$
  By the uniform continuity of $g_{_1},\cdots, g_{_N}$, there exists $\delta>0$ such that for  any two points $x,y$ satisfying $\ud(x,y)<\delta$, we have $$|g_{_i}(x)-g_{_i}(y)|<|\lambda_{n}|\cdot\norm{g_{_i}}_{_{C^0}}, \textrm{for any $i=1,\cdots,N$. }$$


    We choose a sequence of decreasing positive numbers $\{\epsilon_{i}\}_{i\in\N}$ such that $\sum\epsilon_{i}<\delta$ and
    denote $\mathcal{O}_{p_{n}}$ as $\gamma_{n}^{0}$.  Let $\rho$ and $\zeta\in(0,1)$ be the two constants given by Lemma \ref{descend}.
Assume that we already have a periodic orbit $\gamma_{n}^{k}$ which is homoclinically related to $\cO_q$ in $V$,
  then we apply $\gamma_{n}^{k}$ to Lemma \ref{descend}  and we get a periodic orbit $\gamma_{n}^{k+1}$ such that:
  \begin{itemize}
  \item $\lambda^{c}(\gamma_{n}^{k+1})>\zeta\lambda^{c}(\gamma_{n}^{k})$
  \item  $\gamma_{n}^{k+1}$ is a $(\epsilon_{k+1}, 1-\rho|\lambda^c(\gamma_n^k)|)$ good  approximation for $\gamma_{n}^{k}$;
  \item $\gamma_n^{k+1}$ is homoclinically related to $\gamma_n^k$ inside $V$.
  \end{itemize}

  We denote by $\kappa_i=1-\rho|\lambda^c(\gamma_n^i)|$.
   By   induction, we have that   $|\lambda^{c}(\gamma_n^i)|<\zeta^i|\lambda_n|$,
    which implies  that
    $$\kappa_i\geq 1-\rho \cdot\zeta^{i}\cdot|\lambda_n|, \textrm{ for each $i\in \N$}.$$
Therefore, for any integer $k\in\N$,  we have the following estimate:
 $$ 0\geq\sum_{i=0}^{k}\log\kappa_{i}>\sum_{i=0}^{k}2(\kappa_{i}-1)\geq \sum_{i=0}^{k} -2\rho\cdot \zeta^{i}|\lambda_n|>
\frac{2\rho}{1-\zeta}\lambda_{n}.$$
Hence, we have that
 $$ \prod_{i=0}^{k}\kappa_{i}\in(1+\frac{2\rho}{1- \zeta}\lambda_n, 1),$$
 which implies  $\prod_{i=0}^{\infty}\kappa_{i}\in(0,1].$

By Lemma \ref{limit}, $\delta_{\gamma_n^k}$ tends to an ergodic measure $\nu_n$.
Since the center Lyapunov exponent of $\gamma_n^k$ tends to zero when $k$ tends to infinity,
by the continuity of $\log\norm{Df|_{_{E^c}}}$, $\nu_n$ is a non-hyperbolic ergodic measure.

On the other hand, by construction of $\gamma^k_n$, we have that for any $k\in\N$,
 the periodic orbit $\gamma_{n}^{k}$ is a $(\sum_{i=1}^{k}\epsilon_{i},\prod_{i=1}^{k}\kappa_{i})$ good approximation for $\mathcal{O}_{p_{n}}$.
  We denote by $\gamma(n,k)$ the subset of $\gamma_n^k$ corresponding to the one in Definition \ref{good for}.

For any integer $i\in[\,1,N\,]$, we have the following:
 \begin{displaymath} \int g_{_i} \,\ud \delta{\gamma_{n}^{k}}=\frac{1}{\pi(\gamma_{n}^{k})}\sum_{x\in\gamma_{n}^{k}}g_{_i}(x)=
\frac{1}{\pi(\gamma_{n}^{k})\cdot\pi(\gamma_{n}^{0})}
\sum_{x\in\gamma_{n}^{k}}\sum_{j=0}^{\pi(\gamma_{n}^{0})-1}g_{_i}(f^{j}(x)).
\end{displaymath}

 \begin{align*}
 \big|\int g_{_i}\, \ud\delta{\gamma_{n}^{0}}-\int g_{_i}\, \ud\delta{\gamma_{n}^{k}}\big|
 &=\Big|\frac{1}{\pi(\gamma_{n}^{k})\pi(\gamma_{n}^{0})}
\sum_{x\in\gamma_{n}^{k}}\sum_{j=0}^{\pi(\gamma_{n}^{0})-1}g_{_i}(f^{j}(x))
-\frac{1}{\pi(\gamma_{n}^{0})}\sum_{y\in\gamma_{n}^{0}}g_{_i}(y)\Big|
                                                                \\
 &{}=
 \frac{1}{\pi(\gamma_{n}^{k})\pi(\gamma_{n}^{0})}\Big|\sum_{x\in\gamma(n,k)}\big(\sum_{j=0}^{\pi(\gamma_{n}^{0})-1}g_{_i}(f^{j}(x))-
\sum_{y\in\gamma_{n}^{0}}g_{_i}(y)\big)
  \\
&\hspace{5mm}+\sum_{x\in\gamma_{n}^{k}\backslash\gamma(n,k)}\big(\sum_{j=0}^{\pi(\gamma_{n}^{0})-1}
g_{_i}(f^{j}(x))-\sum_{y\in\gamma_{n}^{0}}g_{_i}(y)\big)\Big|
 \\
&{}\leq\frac{1}{\pi(\gamma_{n}^{k})\pi(\gamma_{n}^{0})}\Big( \pi(\gamma_{n}^{k})
\cdot\pi(\gamma_{n}^{0}) \cdot|\lambda_{n}|\cdot  \norm{g_{_i}}_{_{C^0}}
 \\
&\hspace{5mm}+2(1-\prod_{i=1}^{k}\kappa_{i})\cdot\pi(\gamma_{n}^{k})\cdot\pi(\gamma_{n}^{0})\cdot\norm{g_{_i}}_{_{C^0}}\Big)
 \\
&{}\leq|\lambda_{n}|\norm{g_{_i}}+2(1-\prod_{i=1}^{k}\kappa_{i})  \norm{g_{_i}}_{_{C^0}}
\\
&{}\leq(1+\frac{4\rho}{1-\zeta})|\lambda_{n}|\norm{g_{_i}}_{_{C^0}}.
\end{align*}
Hence,  for any $k$, we have that
\begin{align*}
\ud(\delta_{\mathcal{O}_{p_{_n}}},\delta_{\gamma_{n}^{k}})
&= \sum_{i=1}^{\infty}\frac{|\int g_{_i}\,\ud\delta_{\cO_{p_{_n}}}-\int g_{_i}\,\ud\delta_{\gamma_n}|}{2^{i}\norm{g_{_i}}_{_{C^0}}}
\\
&=  \sum_{i=1}^{N}\frac{|\int g_{_i}\,\ud\delta_{\cO_{p_{_n}}}-\int g_{_i}\,\ud\delta_{\gamma_n}|}{2^{i}\norm{g_{_i}}_{_{C^0}}}
+ \sum_{i=N+1}^{\infty}\frac{|\int g_{_i}\,\ud\delta_{\cO_{p_n}}-\int g_{_i}\,\ud\delta_{\gamma_n}|}{2^{i}\norm{g_{_i}}_{_{C^0}}}
\\
&<  \sum_{i=1}^{N}\frac{(1+\frac{4\rho}{1-\zeta})|\lambda_{n}|}{2^{i}}+\sum_{i=N+1}^{\infty}\frac{2}{2^{i}}
\\
&\leq(1+\frac{4\rho}{1-\zeta})|\lambda_{n}|+|\lambda_{n}|.
 \end{align*}

Then, by taking the limit for $k$ tending to infinity,  we get that
$$d(\delta_{\mathcal{O}_{p_{_n}}},\nu_{n})\leq (2+\frac{4\rho}{1-\zeta})\cdot|\lambda_{n}|.$$
Since  $\lambda_{n}$ tends to $0$ and $\delta_{\cO_{p_{_n}}}$ tends to $\mu$,
the non-hyperbolic ergodic measure  $\nu_{n}$ tends to $\mu$.

This ends the proof of Proposition \ref{aa}.
 \endproof
 \begin{Remark}\label{r.aa} From the proof of Proposition \ref{aa},
 one can check  that when $f$ is partially hyperbolic of center dimension one,
 one can  take $V$ to be $M$ and the conclusion of Proposition \ref{aa} still holds.
 \end{Remark} 
\section{Proof of Theorem \ref{thmA} and Proposition \ref{p.non simplex}}
In this section, we first give the  proof of Theorem \ref{thmA}, then using  Theorem \ref{thmA}, we give the proof of Proposition \ref{p.non simplex}.

\subsection{Generation of split flip flop configuration and blender horseshoe}\label{minimality}
In this subsection, we state results on the generation of blender in \cite{BD2}, \cite{BD3}.


\begin{proposition} \label{p.generation of blender}\cite[Proposition 5.6]{BD3} Let $f $be a diffeomorphism with a heterodimensional cycle associated to saddles $P$ and
$Q$ with $ind(P)=ind(Q)+1$. Then there is $g$ arbitrarily $C^1$ close to $f$ exhibiting  a cu-blender horseshoe $\La_g$.
\end{proposition}
From Proposition \ref{p.generation of blender},  we can build blender from co-index one heterodimensional cycle.
To get heterodimensional
cycle, we need the connecting lemma due to  S. Hayashi:

\begin{theorem}\cite{H}
 Let $f\in\diff^1(M)$. For any $C^1$ neighborhood $\mathcal{U}$ of $f$, there is an integer $L=L(\mathcal{U})>0$, such that for any non-periodic
point $z$, there exists two arbitrarily small neighborhoods $B_{z}\subset \tilde{B_{z}}$ of $z$ such that for any two points $x,y\notin\cup_{i=0}^{L}f^{i}(\tilde{B_{z}})$,
if both forward orbit of $x$ and backward orbit of $y$ intersect $B_{z}$, then there exists $g\in\mathcal{U}$ such that $y=g^{n}(x)$, for some integer $n>0$.
\end{theorem}

\subsection{ Assuming minimality of both strong foliations: proof of theorem \ref{thmA}:}
    Now, we can give the proof of theorem A.
    \proof[Proof of Theorem \ref{thmA}]
    By Hayashi's connecting lemma and the transitivity,
    there exists a dense subset of $\mathcal{V}(M)$ such that every diffeomorphism inside this dense subset has a co-index
    one heterodimensional cycle.
    By Propositions ~\ref{p.existence of flip flop} and ~\ref{p.generation of blender}, there exists an open dense subset $\tilde{\cV}(M)$
    of $\cV(M)$ such that for any $f\in\tilde{\cV}(M)$, one has that
     \begin{itemize}
     \item $f$ has a cu-blender horseshoe$(\Lambda^{u}, V^u,\cC^{uu},\mathfrak{D}^u)$ and a cs-blender horseshoe $(\Lambda^{s}, V^s,\cC^{ss},\mathfrak{D}^s)$;
     \item $f$ has a split flip  flop configuration associated to   a dynamically defined $cu$-blender;
     \item $f$ has a split flip flop configuration associated to  a dynamically defined $cs$-blender.
     \end{itemize}
\vspace{2mm}
{\bf Approaching non-hyperbolic ergodic measure by hyperbolic periodic measures}:
Take any hyperbolic ergodic measure $\mu$ and any non-hyperbolic ergodic measure $\nu$,  We will prove that for any $\alpha\in[0,1]$, $\alpha\mu+(1-\alpha)\nu$ is approximated by periodic measures.
By Proposition ~\ref{c1},  the measure $\mu$ is approximated by hyperbolic periodic measures. Hence, we only need to prove it when $\mu$ is a hyperbolic periodic  measure.
Assume that $\mu$ is a periodic measure with positive center Lyapunov exponent (for the negative case, we can argue for the system $f^{-1}$).
Let $\mathcal{O}_{p}$ be the hyperbolic periodic orbit such that  $\mu=\delta_{\cO_p}$.
By the minimality of the strong stable foliation,
$\mathcal{O}_{p}$ is homoclinically related to $\Lambda^{u}$.

      Take a generic point $x$ of measure $\nu$, by the minimality of strong stable foliation, $\mathcal{F}^{ss}(x)$ intersects $W^{u}(\mathcal{O}_{p})$ in  a point $y$. Then, we have that
      \begin{displaymath}
       \lim_{n\rightarrow+\infty}\frac{1}{n}\sum_{j=0}^{n-1}\delta_{f^{j}(y)}=\nu \textrm{ and } \lim_{n\rightarrow+\infty}\frac{1}{n}\sum_{j=0}^{n-1}\norm{Df|_{E^{c}(f^{j}(y))}}=0.
       \end{displaymath}
       By the minimality of strong unstable foliation and the strictly invariant property of $\mathfrak{D}$, for any $\delta>0$, there exists an integer $N_{\delta}$ such that for any strong unstable disc $D^{uu}$ of radius $\delta$, we have that $f^{N_{\delta}}(D)$ contains an element of $\mathfrak{D}^{u}$.  Hence, for any $\delta>0$ and any integer $n\in\N$,
        one has that $f^{N_\delta}(W^{uu}_{\delta}(f^n(y)))$ contains a disc in $\mathfrak{D}$.  Now we can apply the arguments for  Proposition~\ref{p.main} to $\cO_p$ and the point $y$,  proving that
      the invariant measure $\alpha\mu+(1-\alpha)\nu$ is approximated by periodic measures of index $ind(p)$.

       As a consequence, one gets that  every non-hyperbolic ergodic measure is approached by periodic measures of indices $i$ and $i+1$ at the same time.

By the minimality of strong foliations, any two hyperbolic periodic measures of the same index are homoclinically related. As a consequence,
for any two hyperbolic periodic orbits $\gamma_1$ and $\gamma_2$ of the same index, one has that
  $$\overline{\{(1-\alpha)\delta_{\gamma_1}+\alpha\delta_{\gamma_2}: \alpha\in[0,1]\}}\subset \overline{ \cM_{per}(f)}.$$
Hence for any two hyperbolic ergodic measures of same index, by Proposition~\ref{c1},
their convex combination can be approximated by periodic measures of  the same index,
which implies that  $\overline{\mathcal{M}_i(f)}$ and $\overline{\mathcal{M}_{i+1}(f)}$ are convex sets.

   Given two invariant measures $\mu,\nu\in\overline{\mathcal{M}^{*}(f)}$,
   then $\mu$ and $\nu$ are approximated by a sequence of periodic measures
    $\{\delta_{\mathcal{O}_{p_n}}\}$ and $\{\delta_{\mathcal{O}_{q_n}}\}$ of the same index respectively.
    Then the convex combination $\overline{\{(1-\alpha)\delta_{\mathcal{O}_{p_n}}+\alpha\delta_{\mathcal{O}_{q_n}}: \alpha\in[0,1]\}}$
    is contained in the closure of the set of hyperbolic periodic measures.
    Hence, the invariant measure $(1-\alpha)\mu+\alpha\nu$ is accumulated by hyperbolic periodic measures, for any $\alpha\in[\,0,1\,]$.
    On the other hand, for any $\alpha\in[0,1]$, we have $$\int\log\norm{Df|_{E^{c}}}\ud((1-\alpha)\mu+\alpha\nu)=0.$$
     By Proposition \ref{aa} and Remark~\ref{r.aa} , we have that $(1-\alpha)\mu+\alpha\nu$ is accumulated by non-hyperbolic ergodic measure, that is,   $$\{(1-\alpha)\mu+\alpha\nu; \alpha\in[0,1]\}\subset\overline{\cM^{*}(f)}.$$ This proves that $\overline{\cM^{*}(f)}$ is a convex set.

 Since every  non-hyperbolic ergodic measure is  approximated by hyperbolic periodic measures of index $i$ and index $i+1$ at the same time,   we have that
 \begin{itemize}\item $$\overline{\mathcal{M}^{*}(f)}\subset \overline{\mathcal{M}_i(f)}\cap \{\mu\in\mathcal{M}_{inv}(f), \lambda^c(\mu)=0\};$$
 \item $$\overline{\mathcal{M}^{*}(f)}\subset \overline{\mathcal{M}_{i+1}(f)}\cap \{\mu\in\mathcal{M}_{inv}(f), \lambda^c(\mu)=0\};$$
 \item $$\overline{\mathcal{M}^{*}(f)}\subset\overline{\mathcal{M}_i(f)}\cap \overline{\mathcal{M}_{i+1}(f)}.$$
 \end{itemize}

  On the other hand, for any $\mu\in\overline{\mathcal{M}_i(f)}\cap\{\mu\in\mathcal{M}_{inv}(f), \lambda^c(\mu)=0\}$, we have that
   $\mu$ is accumulated by hyperbolic ergodic measures of index $i$.
    By Proposition \ref{c1}, we know that $\mu$ is accumulated by hyperbolic periodic measures of index $i$.
     By Proposition \ref{aa}, the measure $\mu$ is accumulated by non-hyperbolic ergodic measures, ie. $\mu\in\overline{\mathcal{M}^{*}(f)}$.
    Hence, we have
    $$\overline{\mathcal{M}_i(f)}\cap\{\mu\in\mathcal{M}_{inv}(f), \lambda^c(\mu)=0\}\subset \overline{\cM^{*}(f)}.$$
    Similarly, we can prove that
      $$\overline{\mathcal{M}_{i+1}(f)}\cap\{\mu\in\mathcal{M}_{inv}(f), \lambda^c(\mu)=0\}\subset \overline{\cM^{*}(f)}.$$
    Besides, one can  easily  check that
    $$\overline{\mathcal{M}_i(f)}\cap \overline{\mathcal{M}_{i+1}(f)}\subset \overline{\mathcal{M}_i(f)}\cap\{\mu\in\mathcal{M}_{inv}(f), \lambda^c(\mu)=0\};$$
    and
     $$\overline{\mathcal{M}_i(f)}\cap \overline{\mathcal{M}_{i+1}(f)}\subset \overline{\mathcal{M}_{i+1}(f)}\cap\{\mu\in\mathcal{M}_{inv}(f), \lambda^c(\mu)=0\}.$$

    This ends the proof of the  first item.
\vspace{2mm}
    \\
    \textbf{Proof of Item 2}:
   Given $\nu\in\mathcal{M}_{i}(f)$, then there exists a sequence of periodic orbits  $\{\mathcal{O}_{p_n}\}_{n\in\N}$ of s-index $i$ such that $\delta_{\mathcal{O}_{p_n}}$ converges to $\nu$. By the minimality of strong stable foliation and Inclination lemma, for each periodic orbit $\mathcal{O}_{p_n}$, one has that
   \begin{itemize}
   \item the stable manifold of $\cO_{p_n}$ contains a disc in $\cV_{\epsilon_0/2}(\mathfrak{D}^s)$, where $\epsilon_0$ is the strength of $\mathfrak{D}^s$;
   \item there exists a compact submanifold of $\Delta_n$ of $W^u(\cO_{p_n})$ such that each disc in $\mathfrak{D}^s$ intersects $\Delta_n$.
   \end{itemize}
   By Proposition \ref{convex}, there exists an invariant measure $\mu_n$ supported on $\Lambda^{s}$ such that
    for any $\alpha\in[0,1]$, we have that
 $$\alpha\delta_{\mathcal{O}_{p_n}}+(1-\alpha)\mu_n\in\overline{\cM_{per}(f)}.$$
  Let $\xi$ be an accumulation of $\mu_n$, then for any $\alpha\in[0,1]$, we have that  $$\alpha\nu+(1-\alpha)\xi\in \overline{\cM_{per}(f)}.$$

  Similarly, we can prove that for any invariant measure belonging to $\overline{\mathcal{M}_{i+1}(f)}$, there exists an invariant measure supported on $\Lambda^{u}$ such that their convex combination is approximated by periodic measures.

  We only need to take $K_i=\Lambda^{u}$ and $K_{i+1}=\Lambda^{s}$.
 Then the second item is satisfied.

 This ends the proof of Theorem \ref{thmA}.
\endproof
\subsection{$C^1$-generic case: Proof of Proposition \ref{p.non simplex}}

As we know that all the extremal points of $\cM_{inv}(f)$ are ergodic measures.
 In general, a convex subset of $\cM_{inv}(f)$ may have more extreme points.
 Under some assumption, we firstly show that there  exist extreme  points of $\overline{\cM^{*}}(f)$ which are not ergodic.
Recall that $\tilde{\cV}(M)$ is an open and dense subset of $\cV(M)$ given by Theorem~\ref{thmA}.

\begin{lemma}\label{l.extra extremal point}
Given $f\in\tilde{\cV}(M)$. Let $\mu$ and $\nu$ be two hyperbolic ergodic measures of different indices.
Assume that,  for any $\alpha\in[0,1]$, the measure $\alpha\mu+(1-\alpha)\nu$ is approached by hyperbolic periodic measures.

Then there exists $\alpha_{0}\in(0,1)$ such that:
\begin{itemize}\item $$\alpha_0\mu+(1-\alpha_0)\nu\in \overline{\cM^{*}(f)};$$
\item the invariant measure  $\alpha_0\mu+(1-\alpha_0)\nu$ is an extreme point of the convex sets $\overline{\cM^{*}(f)}$,  $\overline{\cM_i(f)}$ and  $\overline{\cM_{i+1}(f)}$.
\end{itemize}
\end{lemma}
\proof Since the indices of $\mu$ and $\nu$ are different,
there exists a unique $\alpha_{0}\in(0,1)$ such that $$\int{\log{\norm{Df|_{E^c}}}}\ud(\alpha_0\mu+(1-\alpha_0)\nu)=0.$$
Since $\alpha_0\mu+(1-\alpha_0)\nu$ is approximated by hyperbolic periodic measures,
by the first  item of Theorem \ref{thmA}, we have that $$\alpha_0\mu+(1-\alpha_0)\nu\in \overline{\cM^{*}(f)}.$$

Assume that  there exist two  measures $\mu_1,\mu_2\in\overline{\cM^{*}(f)}$ and $\beta_0\in(0,1)$ such that
\begin{equation}\label{e.extremal point}
\alpha_0\mu+(1-\alpha_0)\nu=\beta_{0}\mu_1+(1-\beta_0)\mu_2.
\end{equation}

Since $\cM(f)$ is a Choquet simplex, by Equation (\ref{e.extremal point}) and the ergodicity of $\mu$ and $\nu$,
    $\mu_1$ is a convex combination of $\mu$ and $\nu$.
By  the fact that  $\lambda^c(\mu_1)=0$, one has that $\mu_1=\alpha_0\mu+(1-\alpha_0)\nu$,
which implies  that $\mu_1=\mu_2=\alpha_0\mu+(1-\alpha_0)\nu$.
This is proves that   $\alpha_0\mu+(1-\alpha_0)\nu$ is an extreme point of the convex set $\overline{\cM^{*}(f)}$.

   Similarly, we can show that  $\alpha_0\mu+(1-\alpha_0)\nu$ is an extreme point of the convex sets  $\overline{\cM_i(f)}$ and  $\overline{\cM_{i+1}(f)}$.

   This ends the proof of Lemma \ref{l.extra extremal point}.
\endproof
Now, we are ready to give the proof of Proposition \ref{p.non simplex}.
\proof[Proof of Proposition \ref{p.non simplex}]
By Theorem \ref{thmA}, we know that $\overline{\cM^{*}(f)}$, $\overline{\cM_{i}(f)}$ and $\overline{\cM_{i+1}(f)}$ are convex sets.
 By Theorem 3.10  in  \cite{ABC}, there exists a residual subset $\mathcal{R}$ of $\tilde{\cV}(M)$ such that for any $f\in\mathcal{R}$,
 the closure of the set of periodic measures is convex,
which implies that  the convex combinations of hyperbolic ergodic measures of different indices are approached by hyperbolic periodic measures.

We will show that for any $f\in\cR$, none of  $\overline{\cM^{*}(f)}$,  $\overline{\cM_i(f)}$  or  $\overline{\cM_{i+1}(f)}$ is a Choquet simplex.

 We take four hyperbolic ergodic measures $\mu_1,\,\mu_2\in\cM_{i}(f)$ and $\nu_1,\,\nu_2\in\cM_{i+1}(f)$.

  We denote by
$$\cH_1=\{(\alpha_1,\alpha_2,\alpha_3,\alpha_4)\in\mathbb{R}^4| \textrm{ $\alpha_i$ is non-negative and $\Sigma_{i=1}^4 \alpha_i=1$}\}$$
and
 $$\cH_2=\{(\beta_1,\beta_2,\beta_3,\beta_4)\in\mathbb{R}^4| \textrm{ $\beta_{1}\lambda^c(\mu_1)+\beta_2\lambda^c(\mu_2)+\beta_3\lambda^c(\nu_1)+\beta_4\lambda^c(\nu_2)=0$}\}.$$

Since $\lambda^c(\mu_i)>0$ and $\lambda^c(\nu_i)<0$  for $i=1,2$,
one gets that  the hyperplane $\cH_2$ is transverse to $\cH_1$.
As a consequence, the intersection $\cH_1\cap\cH_2$ is a convex quadrilateral whose vertexes are corresponding to four different invariant measures;
Moreover,  each of them is a convex sum of  two hyperbolic ergodic measures of different indices among $\{\mu_1,\mu_2,\nu_1,\nu_2\}$.
 By the generic assumption, these four invariant measures are approximated by hyperbolic periodic measures;
hence by  Theorem \ref{thmA},    they  belong to the set $\overline{\cM^{*}(f)}$.

 By the convexity of the set $\overline{\cM^{*}(f)}$,   the diagonal of  $\cH_1\cap\cH_2$ intersects in  a point which  corresponds to an invariant measure $\mu\in\overline{\cM^{*}(f)}$.  By Lemma \ref{l.extra extremal point},  the vertexes of  $\cH_1\cap\cH_2$ are extreme points of  $\overline{\cM^{*}}(f)$. Hence, $\mu$ is the convex combination of two different pairs of extreme points of $\overline{\cM^{*}(f)}$,
 which implies that $\overline{\cM^{*}(f)}$ is not a Choquet simplex.

  Similarly, one can show that neither $\overline{\cM_i(f)}$ nor  $\overline{\cM_{i+1}(f)}$ is a  Choquet simplex.
\endproof

\section{Invariant measures for Ma\~ n\'e 's  $DA$ example: Proof of theorem \ref{thmB}}\label{s.mane example}

In this section, we first recall Ma\~ n\'e 's  $DA$-example, that is, what exactly the open set $\cW$ in the statement of Theorem~\ref{thmB} is. Then we give the proof the Theorem~\ref{thmB}.
\subsection{Ma\~ n\'e 's $DA$-example}
In \cite{M1}, by doing $DA$ from a linear Anosov diffeomorphism on $\T^3$ whose center is uniformly expanding,
 R. Ma\~ n\'e  constructs an open subset $\cW$ of $\diff^1(\T^3)$ such that for any $f\in\cW$,
the following properties are satisfied:
\begin{itemize}
\item There exist a $Df$-invariant partially hyperbolic splitting $$T\T^3=E^{ss}\oplus E^c\oplus E^{uu}$$
 with $\dim (E^{ss})=\dim (E^c)=1$,
 \item there exists a constant $\lambda\in(0,1)$ such that $\norm{Df^{-1}|_{E^{uu}}}<\lambda;$
\item There exist two hyperbolic periodic orbits of different indices;
\item The center bundle $E^c$ is integrable and the center foliation is minimal. 
\item There exist two open sets $U$ and $V$, a constant $\tau>1$ and  five positive numbers $\epsilon_1,\cdots,\epsilon_5$ such that
\begin{enumerate} \item $\overline{V}$ is a proper subset of $ U$
\item   For any point $x\in\T^3\backslash V$, we have that $$\norm{Df|_{E^{c}(x)}}>\tau;$$
\item Every strong unstable curve of length at least $\epsilon_1$ contains a strong unstable curve of length at least $\epsilon_2$ which is disjoint from $U$;
Moreover, we have that $\lambda\epsilon_1<\epsilon_2$;
    \item Every center plaque of length at least $\epsilon_3$ contains a center plaque of length at least $\epsilon_4$ which is disjoint from $\overline{V}$; Moreover,
    we have that $\tau\epsilon_4>2\epsilon_3$;
        \item For every center leaf $\cF^c(x)$, every connected component of $\cF^c(x)\cap(U\backslash\overline{V})$ has length larger than $\epsilon_5$; Moreover, we have that  $\tau\epsilon_5>\epsilon_3$.
\end{enumerate}
\end{itemize}
By construction, the diffeomorphism $f$ is isotopic to  a linear Anosov.

R. Ma\~ n\'e   proved the followings:
\begin{theorem}\label{thm.mane}\cite[Theorem B]{M1} For every $f\in\cW$, the diffeomorphism $f$ is robustly transitive and non-hyperbolic.
\end{theorem}

\begin{lemma}\label{M}\cite[Lemma 5.2]{M1}
   Let $f\in\cW$. For any  $x\in\mathbb{T}^{3}$, there exists $y\in W^{uu}_{\epsilon_1}(x)$ such that the forward orbit of $y$ is contained in $\mathbb{T}^{3}\backslash U$.
 \end{lemma}
A recent result by R.Potrie ~\cite{Po} implies  that $f\in\cW$ is dynamically coherent, that is, there exist invariant foliations $\cF^{cs}$ and $\cF^{cu}$ tangent to $E^{ss}\oplus E^c$ and $E^c\oplus E^{uu}$ respectively.

\subsection{Proof of Theorem \ref{thmB}}
We deal with Ma\~ n\'e 's example separately because we only know that strong stable  foliation and center foliation are minimal. The minimality of strong unstable foliation is still unknown.

The minimality of strong stable foliation is due to \cite{BDU,PS}:

 \begin{lemma}\label{mini} \cite{BDU,PS} There exists an open dense subset $\cW^{s}$ of $\cW$ such that for any $f\in\cW^{s}$, the strong stable foliation is minimal.
 \end{lemma}
 The proof  of theorem B strongly depends on the properties of $\cW$, ie. the DA construction of  Ma\~ n\'e 's  example.

 \proof[Proof of Theorem \ref{thmB}]
  By Hayashi's connecting lemma and transitivity, there exists  a dense subset of $\mathcal{W}$ such that every diffeomorphism in  this dense subset has a co-index one heterodimensional cycle.
 By Propositions   ~\ref{p.existence of flip flop} and \ref{p.generation of blender},
there exists an open dense subset of $\cW$ such that every diffeomorphism in this set has a split flip flop configuration associated to a dynamically defined $cs$-blender.
  On the other hand, by Lemma \ref{mini}, there exists an open dense subset  $\cW^{s}$  of $\mathcal{W}$
such that the strong stable foliation is minimal for any $f\in\cW^s$.

    To sum up, there exists an open and dense subset $\tilde{\cW}$ of $\cW$ such that for any $f\in\tilde{\cW}$, we have that:
 \begin{itemize}
  \item $f$ has minimal strong stable foliation;
 \item $f$ has a split flip flop configuration associated to a dynamically defined $cs$-blender.
 \end{itemize}
Now, we fix a $Df$ strictly invariant strong unstable cone field $\cC^{uu}$ around $E^{uu}$ on $\T^3$ 
such that $Df$ is uniformly expanding along $\cC^{uu}$.

\textbf{Non-hyperbolic ergodic measures approached by hyperbolic periodic measures of index one}.

We will prove that for any hyperbolic periodic measure of index $1$ and any non-hyperbolic ergodic measure, their convex combination is approximated by periodic measures.

By  construction, there exists $\delta_0>0$ such that $B_{\delta_0}(V)\subset U$.
 We fix a hyperbolic periodic point $p$ of index $1$, the size of local unstable manifold $W^{u}_{loc}(\cO_p)$   and the size of local strong stable manifold $W^{s}_{loc}(\cO_p)$.
 By the minimality of the  foliation $\mathcal{F}^{ss}$ and the fact that  $W^s(\cO_p)=W^{ss}(\cO_p)$,
there exists a positive  integer $k$ such that
 \begin{itemize}
 \item  the integer $k$ only depends on the size  of $W^{s}_{loc}(p)$ and $\delta_0$;
 \item  for any $x\in\mathbb{T}^{3}$, $f^{-k}(W^{s}_{loc}(p))$ intersects $\cF^{cu}_{\delta_0}(x)$ transversely,  where $\cF^{cu}_{\delta_0}(x)$  denotes the $\delta_0$ neighborhood of $x$ in the leaf $\cF^{cu}(x)$.
     \end{itemize}

 Given a non-hyperbolic ergodic measure $\nu\in\mathcal{M}_{erg}(\mathbb{T}^{3},f)$, and consider the convex sum $\alpha\delta_{\mathcal{O}_{p}}+(1-\alpha)\nu$, for $\alpha\in[0,1]$.
Then  for any $\alpha\in(\,0,1\,]$, we have that
$$\lambda^{c}(\alpha\delta_{\mathcal{O}_{p}}+(1-\alpha)\nu)>0.$$
    We fix $\alpha\in(0,1]$ and we denote by  $$\lambda^{\prime}=\exp(-\lambda^{c}(\alpha\delta_{\mathcal{O}_{p}}+(1-\alpha)\nu)).$$
     Lemma \ref{shadow} provides two positive numbers $L$ and $d_{0}$,  corresponding to the number $\sqrt{\lambda^{\prime}}$ and to  the splitting $E^{ss}\oplus (E^c\oplus E^{uu})$, such that for any $d\in(0, d_0)$, we have that any $\sqrt{\lambda^{\prime}}$-quasi hyperbolic periodic $d$-pseudo orbit is $L\cdot d$ shadowed by a periodic orbit.

     We choose a number $d\in(0,d_{0})$ such that $L\cdot d$ is small enough, whose precise value would be fixed at the end.
 Then there exists an integer $N_{d}$ such that
$$f^{-N_{d}}(W^{u}_{loc}(\mathcal{O}_{p}))\subset W^u_{d/2}(\mathcal{O}_{p}) \textrm{ and  } f^{N_{d}}(W^{s}_{loc}(\mathcal{O}_{p}))\subset W^s_{d/2}(\mathcal{O}_{p}).$$

Now, we fix a point $x$ in the basin of $\nu$; Since  strong stable foliation is minimal and $W^{u}_{loc}(p)$ is everywhere tangent to $E^c\oplus E^{uu}$,  there exists a transversely intersection $y$ between the strong stable manifold of $x$ and  $W^{u}_{loc}(\mathcal{O}_{p})$.
By the choice of $y$,  we have that
$$\lim_{n\rightarrow+\infty}\frac{1}{n}\sum_{i=0}^{n-1}\delta_{f^{i}(y)}=\nu \textrm{ and } \lim_{n\rightarrow+\infty}\frac{1}{n}\sum_{i=0}^{n-1}\log\norm{Df|_{E^{c}(f^{i}(y))}}=0.$$

In the following,  we will find a $\sqrt{\lambda^{\prime}}$-quasi hyperbolic periodic $d$-pseudo orbit, with large period,  such that it  spends almost $1-\alpha$ proportion of its period to  follow  the forward orbit of $y$  and also spends almost $\alpha$ proportion of its period to  follow the periodic orbit $\cO_p$.

      Take $\epsilon<\frac{-\log{\lambda^{\prime}}}{16}$ small, then we have the following:
      \begin{itemize} \item there exists $\delta>0$ such that for any two points $z,w\in\mathbb{T}^3$ satisfying $\ud(z,w)<\delta$, we have that
       $$ -\epsilon \leqslant \log{\norm{Df|_{{E}^{c}(w)}}}-\log{\norm{Df|_{{E}^{c}(z)}}}\leqslant  \epsilon ;$$
      \item there exists  an integer $N$ such that for any $n\geqslant N$, we have the following:
      \begin{displaymath}
      -\epsilon< \frac{1}{n}\sum_{i=0}^{n-1}\log\norm{Df|_{E^{c}_{f^{i}(y)}}}<\epsilon \textrm{ and  } \ud(\,\frac{1}{n}\sum_{i=0}^{n-1}\delta_{f^{i}(y)},\,\nu)<\epsilon/2.
      \end{displaymath}
      \end{itemize}

For any $C^1$ curve $\gamma$ in $\mathbb{T}^3$, we denote by $\ell(\gamma)$ the length of $\gamma$.

 We choose a $C^1$ curve $\gamma_{n}\subset W^{c}(y)$, centered at $y$,  such that  $\ell(\gamma_{n})=\delta e^{-2n\epsilon}$. Now, consider the set  $S_{n}$ which is the  $\ell(\gamma_n)$ tubular neighborhood of  $W^{uu}_{\delta}(x)$ in the leaf $\cF^{cu}(x)$.

   Similar to Claim~\ref{c.length of central curve}, for  $n$  large enough,  one has that
 $$\ell(f^{i}(\gamma_{n}))\leq \norm{Df}^{i}\ell(\gamma_{n})<\delta,\textrm{ for any $i=0,\cdots n$}.$$

    According to this estimate and the choice of $\delta$, we have that $$\ell(f^{n}(\gamma_{n}))=\int_{0}^{1}\|\frac{\ud}{\ud t}f^{n}(\gamma_{n}(t))\|\ud t
    \geq \int_{0}^{1}e^{-n\epsilon}\prod_{j=0}^{n-1}\norm{Df|_{E^{c}(f^{j}(y))}}\|\gamma_{n}^{\prime}(t)\|\ud t\geq \delta e^{-4n\epsilon}.$$

    Consider the connected component of $f^{n}(S_{n})\cap B_{\delta}(f^{n}(y))$ which contains $f^{n}(y)$, and we  denote it as $\tilde{S_{n}}$;   then $\tilde{S_{n}}$ is a $C^1$ surface tangent to $E^c\oplus E^{u}$ satisfying that
    \begin{itemize}
    \item $f^n(\gamma_n)\subset \tilde{S_n}$;
    \item $\tilde{S_n}$ is foliated by curves  tangent to $\cC^{uu}$, whose lengths are  $2\delta$ ;
    \item For any center plaque $\gamma\subset\tilde{S_n}$  joining the two  boundary components of $\tilde{S_n}$ which are tangent to $\cC^{uu}$,
    one has that
     $\ell(\gamma)\geq\delta e^{-5n\epsilon}$.
     \end{itemize}
     Now, we will iterate $\tilde{S_{n}}$ to make it cut the local stable manifold of $\mathcal{O}_{p}$.
By the uniform expansion in the strong unstable direction,
 there exists an integer $N_{\delta}$ such that for any curve $D^{uu}$ tangent to $\cC^{uu}$ of length at least $\delta/2$,
  we have that
                                               $$\ell(f^{N_{\delta}}(D^{uu}))\geq \epsilon_1.$$
     Hence, $f^{N_{\delta}}(\tilde S_{n})$ is foliated by curves tangent to $\cC^{uu}$,  whose  lengths are  no less than $4\epsilon_1$.
     By Lemma \ref{M}, there exists $z\in W^{uu}_{\delta/2}(f^{n}(y))$ such that $Orb^{+}(f^{N_{\delta}}(z))$ is contained in $\mathbb{T}^3\backslash B_{\delta_0}(V)$.

      We denote by $W_n$ the connected component of  $f^{N_{\delta}}(\tilde S_{n})\cap B_{\delta/2}(f^{N_{\delta}}(z))$ containing $f^{N_{\delta}}(z)$.  Let $[\cdot]$ denote the integer part of a number, we denote by
      $$T_{n}=\Big[\frac{5\epsilon n+\log\delta_0-\log\delta+N_{\delta}\cdot\log{b}}{\log\tau}\Big]+1,\textrm{ where $b>\sup_{x\in\T^3}\norm{Df^{-1}(x)}$};$$
      Since for any point $x\in\mathbb{T}^3\backslash V$, we have
      $$\norm{Df|_{_{E^c(x)}}}\geq \tau,$$
      hence $f^{T_{n}}(W_{n})$ contains a disc tangent $E^c\oplus E^{u}$ whose diameter is no less than $\delta_{0}$.
      By the choice of $k$,   $f^{T_{n}}(W_{n})$  intersects $f^{-k}(W^{s}_{loc}(p))$ transversely.
      We denote by $$t_n=2N_{d}+k+N_{\delta}+T_n.$$
      To sum up,  there exists a point $w\in W^{u}_{d/2}(p)$ such that
      \begin{itemize}
      \item  $$f^{n+t_n}(w)\in W^{s}_{d/2}(p);$$
      \item the orbit segment $\{f^{N_d}(w),\cdots,f^{n+N_d}(w)\}$ follows the orbit segment $\{y,\cdots, f^{n}(y)\}$ at a distance less than $\delta$.
      \end{itemize}
      When we choose $n$ large, we have that $$\frac{t_n}{n}\leq\frac{6\epsilon}{\log\tau}.$$
      Moreover, we have the following:
     $$\big|\frac{1}{n+t_n}\sum_{i=0}^{n+t_n-1}\log\norm{Df|_{E^{c}_{f^{i}(w)}}}\big|<\epsilon,$$
    $$\ud(\frac{1}{n+t_n}\sum_{i=0}^{n+t_n-1}\delta_{f^{i}(w)},\:\nu) <\frac{6}{\log\tau}\,\epsilon+\epsilon.$$

\begin{claim} There exist integers $n$ and $m$ arbitrarily  large such that
\begin{itemize}
\item $$\big|\frac{m\pi(p)}{n}-\frac{\alpha}{1-\alpha}\big|+\frac{2N_d}{n}<\epsilon;$$
\item the orbit segment $\{w,t_n+n+m\pi(p)\}$ is a $\sqrt{\lambda^{\prime}}$-quasi hyperbolic string corresponding to the splitting $E^{ss}\oplus (E^c\oplus E^{uu})$.
\end{itemize}
\end{claim}

      \begin{figure}[h]
\begin{center}
\def\svgwidth{0.9\columnwidth}
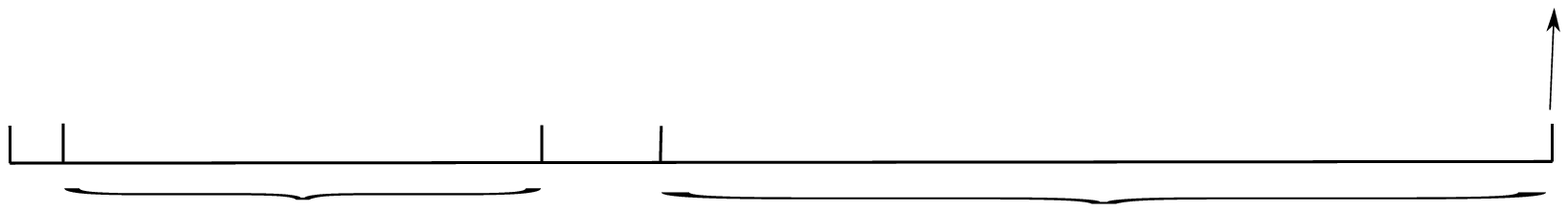
  \caption{}
\end{center}
\end{figure}
The proof of the claim above is similar to the one of   Claim~\ref{c.choice of  t_k}.
  Once again,  by Lemma \ref{shadow}, we have a periodic orbit $\cO_{p^{\prime}}$ of index one which shadows the orbit segment  $\{w,t_n+n+m\pi(p)\}$  in a distance $L\cdot d$. Moreover, when $d$ is chosen small, one has that
   $$\ud(\delta_{\cO_{p^{\prime}}},\alpha\delta_{\mathcal{O}_{p}}+(1-\alpha)\nu)<c\cdot\epsilon, \textrm{where $c$ is a constant independent of $\epsilon$}. $$

    By the arbitrary choice of $\alpha$ and
    compactness of the set $\{\alpha\delta_{\mathcal{O}_{p}}+(1-\alpha)\nu|\alpha\in[\,0,1\,]\}$, the ergodic measure $\nu$ is  approximated
    by hyperbolic ergodic measures of index $1$.
    \\
    \textbf{Convexity of the set $\overline{\cM^{*}(f)}$}.
    Since any non-hyperbolic ergodic measure is approximated by hyperbolic periodic measures of index $1$,
    we have that for any $\mu,\nu\in\overline{\cM^{*}(f)}$, both $\mu$ and  $\nu$ are approximated by hyperbolic periodic measures of index $1$.
     Since the hyperbolic periodic orbits of index $1$ are homoclinically related to each other,  one has that $$\{\alpha\mu+(1-\alpha)\nu;\alpha\in[0,1]\}$$
     is contained in the closure of the set of hyperbolic periodic measures of index $1$.
     Notice that $\alpha\mu+(1-\alpha)\nu$ has  zero mean  center Lyapunov exponent,
      for any $\alpha\in[0,1]$.
     By Proposition \ref{aa} and Remark \ref{r.aa}, we have that $\alpha\mu+(1-\alpha)\nu$ is approximated by non-hyperbolic ergodic measure. Hence, we have that
     $$\{\alpha\mu+(1-\alpha)\nu;\alpha\in[0,1]\}\subset\overline{\cM^{*}(f)}.$$

This ends the proof of Theorem ~\ref{thmB}.
    \endproof

{\bf Acknowledgment: }
We would like to  thank J\'er\^ome Buzzi, Sylvain Crovisier, Shaobo Gan, Lan Wen and Dawei Yang   for useful comments.

Jinhua Zhang  would like to thank Institut de Math\'ematiques de
Bourgogne for hospitality and China Scholarship Council (CSC) for financial support (201406010010).

\bibliographystyle{plain}

\vskip 2mm

\noindent Christian Bonatti,

\noindent {\small Institut de Math\'ematiques de Bourgogne\\
UMR 5584 du CNRS}

\noindent {\small Universit\'e de Bourgogne, 21004 Dijon, FRANCE}

\noindent {\footnotesize{E-mail : bonatti@u-bourgogne.fr}}

\vskip 2mm
\noindent Jinhua Zhang,

\noindent{\small School of Mathematical Sciences\\
Peking University, Beijing 100871, China}

\noindent {\footnotesize{E-mail :zjh200889@gmail.com }}\\
\noindent and\\
\noindent {\small Institut de Math\'ematiques de Bourgogne\\
UMR 5584 du CNRS}

\noindent {\small Universit\'e de Bourgogne, 21004 Dijon, FRANCE}\\
\noindent {\footnotesize{E-mail : jinhua.zhang@u-bourgogne.fr}}


\end{document}